\documentclass{cta-author}
\usepackage{epstopdf}
\usepackage{amsmath}
\usepackage{cases}
\usepackage[caption=false,font=footnotesize]{subfig}
\usepackage{fixltx2e}
\usepackage{amssymb}
\usepackage{mathtools}
\usepackage{bm}
\usepackage{multirow}
\usepackage{color}
\usepackage{algorithm}
\usepackage{algorithmic}
\usepackage{amsthm}
\usepackage{mathrsfs}
\usepackage{enumitem}
\usepackage[euler]{textgreek}
\usepackage{subfig}

{}
{}
{}

\begin{document}

\supertitle{Accepted by IET Renewable Power Generation}

\title{Tertiary Regulation of Cascaded Run-of-the-River Hydropower in the Islanded Renewable Power System Considering Multi-Timescale Dynamics}

\author{\au{Yiwei Qiu$^{1,2}$} \au{Jin Lin$^{1\corr}$} \au{Feng Liu$^{1}$} \au{Ningyi Dai$^{2}$} \au{Yonghua Song$^{2,1}$} \au{Gang Chen$^{3}$} \au{Lijie Ding$^{3}$}}

\address{\add{1}{State Key Laboratory of Control and Simulation of Power Systems and Generation Equipment, Department of Electrical Engineering, Tsinghua University, Beijing 100087, China}
\add{2}{State Key Laboratory of Internet of Things for Smart City, University of Macau, Macau 999078, China}
\add{3}{State Grid Sichuan Electric Power Research Institute, Chengdu 610000, China}
\email{linjin@tsinghua.edu.cn}}

\begin{abstract}
\looseness=-1
  To enable power supply in rural areas and to exploit clean energy, fully renewable power systems consisting of cascaded run-of-the-river hydropower and volatile energies such as pv and wind are built around the world. In islanded operation mode, the primary and secondary frequency control, i.e., hydro governors and automatic generation control (AGC), are responsible for the frequency stability. However, due to limited water storage capacity of run-of-the-river hydropower and river dynamics constraints, without coordination between the cascaded plants, the traditional AGC with fixed participation factors cannot fully exploit the adjustability of cascaded hydropower. When imbalances between the volatile energy and load occur, load shedding can be inevitable. To address this issue, this paper proposes a coordinated tertiary control approach by jointly considering power system dynamics and the river dynamics that couples the cascaded hydropower plants. The timescales of the power system and river dynamics are very different. To unify the multi-timescale dynamics to establish a model predictive controller that coordinates the cascaded plants, the relation between AGC parameters and turbine discharge over a time interval is approximated by a data-based second-order polynomial surrogate model. The cascaded plants are coordinated by optimising AGC participation factors in a receding-horizon manner, and load shedding is minimised. Simulation of a real-life system with real-time pv data collected on site shows the proposed method significantly reduces load loss under pv volatility.
\end{abstract}


\maketitle

\section*{Nomenclature}

\begin{description}[leftmargin=4em,style=nextline]

\setlength{\parskip}{1.5pt}

\item[$\tilde{x}$, $\hat{x}$] The nominal/initial value and the increment of variable $x$, i.e., $x = \tilde{x} + \hat{x}$
\item[$\underline{x}$, $\overline{x}$] Lower and upper limits of $x$
\item[$\bm{x}$] Vector with entries $x_i$
\item[$\omega^{\mathrm{ref}}$] System frequency deviation reference
\item[$\Omega^{\mathrm{G}}$] Set of hydropower generators
\item[$\Omega^{\mathrm{sh}}$] Set of sheddable loads

\item[$P_{\mathrm{G}i}$] Power of the $i$th hydropower unit
\item[$P_{\mathrm{S}}$] Power output of volatile energies

{\color{black}

\item[$P_{\mathrm{S}}^{\mathrm{avail}}$] Available volatile energies
\item[$P_{\mathrm{S}}^{\mathrm{pred}}$] Prediction of volatile energies
\item[$P_{\mathrm{S}}^{\mathrm{cur}}$] Curtailment of volatile energies

}

\item[$P_{\mathrm{D}}$] Load power
\item[$c_{\mathrm{G}i}$] Participation factor of the $i$th generator in AGC
\item[$P_\mathrm{D}^\mathrm{sh}$] Total load shedding
\item[$P_{\mathrm{D}j}^\mathrm{sh}$] Power of the $i$th sheddable load
\item[$\pi_i$] Binary variable of the $i$th sheddable load
\item[$\bm{p}$] $\triangleq [\bm{c}_{\mathrm{G}}^\mathrm{T}, \omega^{\mathrm{ref}},P_\mathrm{D}^\mathrm{sh} ]^\mathrm{T}$, vector of tertiary control variables

\item[$\bm{P}_{\mathrm{Inj}}$] Nodal power injection vector
\item[$\bm{P}_\mathrm{B}$] Vector of the branch power flow
\item[$\bm{A}_\mathrm{Net}$] Susceptance-weighted incidence matrix
\item[$\bm{B}_\mathrm{Net}$] Nodal susceptance matrix

\item[$Q^{\mathrm{turb}}_{\mathrm{H}i}$] Turbine discharge at the $i$th hydro plant
\item[$Q_{\mathrm{H}i}^\mathrm{sp}$] Water spillage at the $i$th hydro plant
\item[$H_{\mathrm{H}i}^\mathrm{up}$] Upstream water stage at the $i$th hydro plant
\item[$H_{\mathrm{H}i}^\mathrm{down}$] Downstream water stage at the $i$th hydro plant
\item[$\eta_i$] Efficiency of the $i$th hydropower unit
\item[$g$] Gravity acceleration
\item[$H$, $Q$] Water stage and discharge along the river
\item[$I_0$, $I_f$ ] River bed and equivalent friction slopes.
\item[$\Omega^\mathrm{m}$] Set of monitoring points on the river

\item[$\bm{x}_\mathrm{R}$] $\triangleq \big[ \hat{\bm{Q}}^\mathrm{T}, \hat{\bm{H}}^\mathrm{T} \big] \hspace{0pt}^\mathrm{T}$, river state vector.
\item[$\bm{u}_\mathrm{R}$] $\triangleq \big[  \bm{Q}_{\mathrm{H}}^\mathrm{turb} \hspace{0pt}^\mathrm{T}, \bm{Q}_{\mathrm{H}}^\mathrm{sp}\hspace{0pt}^\mathrm{T} \big] \hspace{0pt}^\mathrm{T}$, river control vector
\item[$\bm{z}_\mathrm{R}$] River boundary condition vector, including upstream inflow and downstream stage
\item[$\phi_i(\cdot)$] Legendre polynomial

\end{description}

\section{Introduction}
\label{sec:intro}

\begin{figure}[tb]
	\centering
	\includegraphics[width=3.15in]{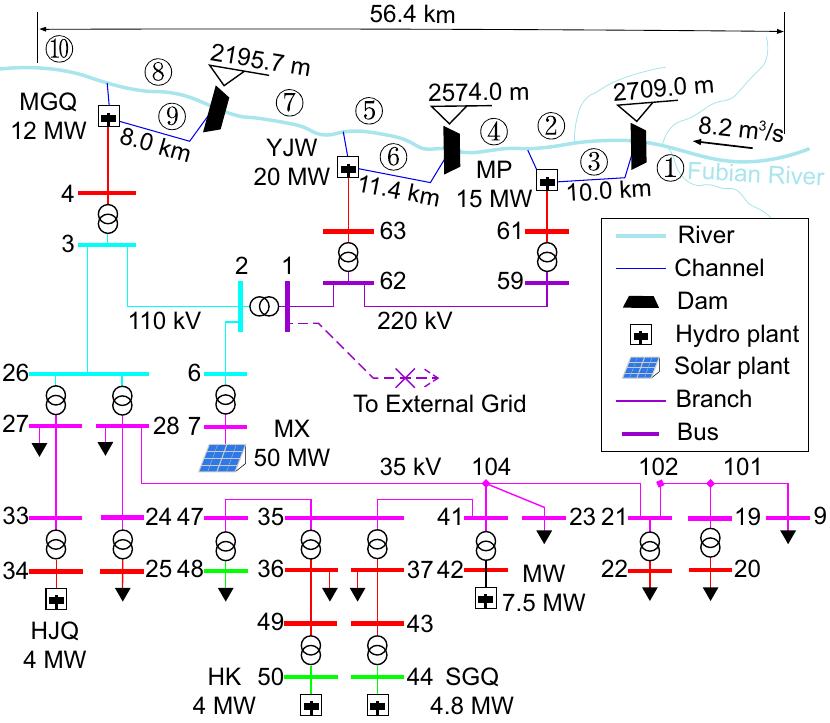}
	\caption{Diagram of a typical power system with three cascaded run-of-the-river hydropower plants and a pv plant in islanded operation mode, which is located in Xiaojin County, Sichuan Province, China}
	\label{fig:xiaojin}
\end{figure}

In  mountain areas such as in Tibet and Sichuan Province in China as well as the plateaus in South Asia and Africa, to enable electrical power supply to the local residents and to exploit clean energy, cascaded run-of-the-river hydropower plants have been built along river valleys \cite{chen2020emergy,wang2019research,bhandari2011electrification,kasman2019performance,sterl2018new}. Moreover, in recent years, photovoltaic (pv) and wind plants are built in these areas to further exploit the renewable energy and to compensate for power shortages in dry seasons \cite{yang2010research,Jiang2015Growt,wang2019research}. The hydroelectric, solar, and wind power make up fully renewable local power systems. A typical example located in Xiaojin County, Sichuan Province, China, is shown in Fig. \ref{fig:xiaojin}.

Among these power systems, some are designed for islanded operation \cite{kaldellis2001evaluation,papaefthymiou2010wind}, and some are connected to the external grid via long-distance transmission lines in normal operation mode, but these need to operate in islanded mode during planned or accidental transmission line outages, especially those caused by natural disasters in rural mountain areas \cite{wang2019research,bhandari2011electrification}. For instance, the Xiaojin system shown in Fig. \ref{fig:xiaojin} was forced to operate in islanded mode in June, 2020 due to that a landslide cut off the transmission line to the external grid.

In islanded operation mode, considering the volatility of solar, wind, and load power, energy balance and frequency stability are the most important issues \cite{martinez2018frequency,martinez2016frequency}.
In this situation,
due to the lack of support from an external power grid and energy storage systems (ESSs) in rural areas, the primary and secondary control of cascaded hydropower, i.e., the hydro turbine governors and automatic generation control (AGC), are essential to stabilise the frequency against pv/wind and load volatilities. 
{\color{black}
But with a high penetration of renewables, frequency stability is challenging, as the pv and wind volatility may exceed the adjustability of hydropower units which may cause stability issues.
Although by simply disconnecting pv or wind from the islanded grid instability can be completely avoided, it may cause power supply shortage and consequent loss of load. Considering socioeconomic benefits, applying controlled curtailment for smoothing renewable power \cite{Howlader2015Integrated,howlader2020active,Ma2019optimal} and using the hydropower to regulate the system seems to be a better option.}

Until now, on the topic of using cascaded hydropower to mitigate solar and wind volatility, elaborations have been made by the community, including scheduling \cite{Jaramillo2004,Apostolopoulou2018Optimal,Apostolopoulou2018robust} and online control \cite{Hug2011Predictive, Hamann2016Using,qiu2020stochastic}, but these works focus on grid-connected operations. On the other hand, studies on the frequency control of islanded power systems with hydropower focus on the timescale of the electromechanical dynamics, relating to the primary and secondary control \cite{martinez2018frequency,martinez2016frequency}.

In addition to the primary and secondary frequency control, in a longer time scale the mismatch between the stochastic pv power and load accumulates. If the accumulated energy is not properly allocated to the cascaded hydropower plants, due to river dynamics some of them may lose adjusting ability as explained later, endangering stability or causing loss of load. To ensure normal operation over a longer timescale, a study on the tertiary regulation that considers the hydraulic coupling between the cascaded plants is also needed, which is, to the best of the authors' knowledge, still lacking. This paper aims to fill this gap.

In stabilising the frequency against the volatilities of wind and solar energy, the primary and secondary frequency control adjusts the power output of the cascaded hydropower plants.
Forthermore, in addition to the capacity and ramping limits, the adjustment of cascaded run-of-the-river hydropower is also limited by river dynamics. This is because in contrast to the conventional dam hydropower, in run-of-the-river hydropower, the water storage capacity is very limited, and the water energy is spatially distributed along the river. The utilisation of water is subject to the river dynamics, and the cascaded plants are therefore hydraulically coupled \cite{Hug2011Predictive, Hamann2016Using,qiu2020stochastic}. Moreover, river operation and ecological regulations often require that the water stage along the river is limited within an allowed interval \cite{Link2010model}, which further limits the adjustment of hydropower.

If we do not consider the river dynamics but only use a traditional AGC that allocates incremental power with fixed proportions to the cascaded plants, due to the impact of solar and wind volatility on hydropower generation and then the water distribution along the river, unacceptable violation of the river operation constraints or a large amount of load shedding may occur, as exemplified in Section \ref{sec:benchmark}. In contrast, if we consider the river dynamics and accordingly coordinate the cascaded plants in a tertiary regulation scheme by dynamically adjusting the AGC participation factors, as later illustrated in Fig. \ref{fig:frame}, the adjustment ability of the hydropower can be maximised, and consequently load loss can be minimised.

\begin{figure}[tb]
   {\centering
    \subfloat[]{\includegraphics[width=3.3in]{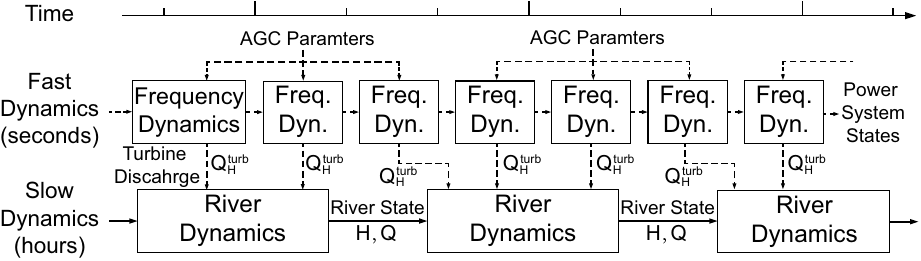}} \\
    \subfloat[]{\includegraphics[width=3.3in]{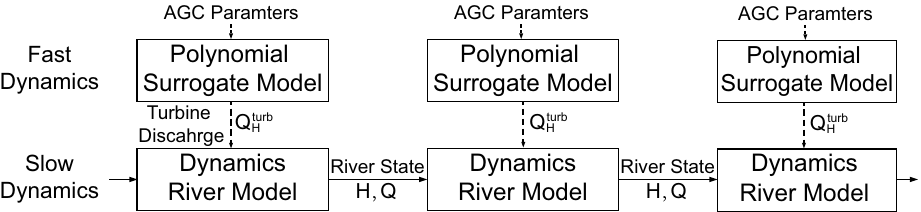}} \\ }
  \caption{Illustration of deling with the multi-timescale dynamics of the power system and cascaded run-of-the-river hydropower}
    \figfooter{a}{The multi-timescale dynamics of the actual system}
    \figfooter{b}{Modeling of the dynamics in the proposed tertiary regulation method}
	\label{fig:timescale}
\end{figure}

\begin{figure*}[tb]
	\centering
	\includegraphics[width=5.38in]{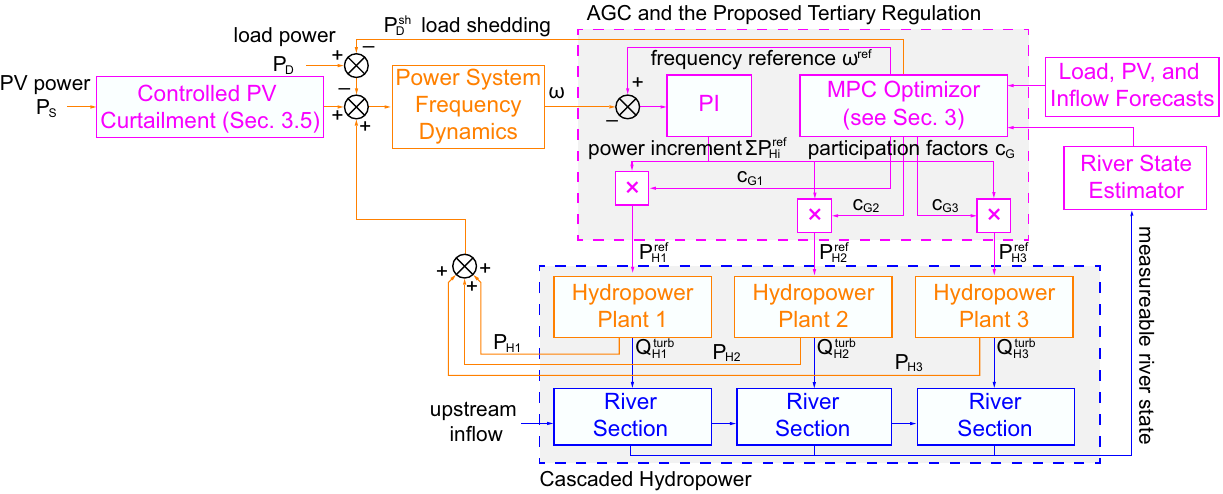}
	\caption{Framework of the proposed regulation approach. The orange, blue, and purple colors represent the electromechanical side, the hydraulic side, and the controller, respectively}
	\label{fig:frame}
\end{figure*}

However, power system frequency dynamics and river dynamics, which is usually characterised by the \emph{shallow water equations} \cite{Link2010model,Litrico2009,Hug2011Predictive,Hamann2016Using,qiu2020stochastic}, have very different timescales, as illustrated in Fig. \ref{fig:timescale}(a).
If we directly combine them to establish a model predictive controller (MPC) for the tertiary regulation with a time resolution compatible with the frequency dynamics and a horizon that can accommodate the river dynamics, the \emph{curse of dimensionality} will arise. To address this issue, in this work a data-based polynomial surrogate model describing the relation between the tertiary control variables (including AGC participation factors, system frequency reference and load shedding) and turbine discharges over a certain time interval is developed.
Thus, the detailed model of power system dynamics can be replaced by simple algebraic functions and then easily incorporated into the MPC formulation with the dynamic river model without causing the curse of dimensionality, as illustrated in Fig. \ref{fig:timescale}(b).

Following the above presented ideas, this paper proposes a coordinated tertiary regulation approach for islanded power systems with cascaded run-of-the-river hydropower and volatile generations.
The framework of the proposed approach is shown in Fig. \ref{fig:frame}. Specifically, the following two contributions are made:
\begin{enumerate} [leftmargin=2em,style=nextline]
  \item A data-based polynomial surrogate model is first established to describe the hydro turbine discharge as a function of AGC participation factors, power system frequency reference, and load shedding. This model can be easily incorporated into an MPC formulation for the tertiary regulation of cascaded run-of-the-river hydropower.
  \item Based on the polynomial surrogate models of power system dynamics and a dynamic river model based on the shallow water equations, a tertiary regulation approach based on MPC jointly considering the multi-timescale dynamics is proposed by updating the AGC participation factors to coordinate the cascaded plants and to minimise load shedding.
\end{enumerate}

Simulation of a real-life system with actual real-time pv data verifies that compared to other regulation schemes the proposed approach significantly reduces load loss.

Specially note that our previous work \cite{qiu2020stochastic} focuses on grid-connected operation of the cascaded hydro-solar system and does not consider any multi-timescale dynamics. Despite both works adopt the same river model, the scope, methodology, and contribution of \cite{qiu2020stochastic} and this work are totally different.

This paper is organised as follows. Section \ref{sec:model} develops the polynomial surrogate model for the electromechanical side of the islanded system and introduces the dynamic river model. Section \ref{sec:controller} presents the optimal control formulation of the tertiary regulation. Finally, in Section \ref{sec:case}, case studies are presented.

\section{Dynamic Modeling of Islanded Power Systems with Cascaded Run-of-the-River Hydropower}
\label{sec:model}

To coordinate cascaded hydropower plants in the tertiary regulation, the electromechanical side, i.e., the power system dynamics, and the hydraulic side, i.e., the river dynamics, need to be jointly considered. As noted in the Introduction, we establish a polynomial surrogate model in Section \ref{sec:modelelectric} to describe the relation between the water discharge of the hydropower plants and the tertiary control variables, which represents the impact of the primary and secondary frequency control on the hydraulic side. Section \ref{sec:river} uses the shallow water equations to depict the river dynamics. Then, in Section \ref{sec:overallmodel}, these models are combined to depict the overall system.

\subsection{Electromechanical Side: Polynomial Surrogate Model of Power System Dynamics}
\label{sec:modelelectric}

{\color{black}With proper smoothing curtailment control of renewable generations, the frequency of an islanded power system with small hydropower and volatile energy can be stabilised by turbine governors and AGC \cite{
kaldellis2001evaluation,papaefthymiou2010wind,martinez2018frequency,martinez2016frequency}, known as the \emph{primary} and \emph{secondary control}.}
However, considering the river dynamic constraints in the cascaded run-of-the-river hydropower, the AGC arguments, for instance the participation factors, should be repeatedly optimised to coordinate the cascaded plants over longer time periods, which is referred to as the \emph{tertiary regulation} or \emph{tertiary control} in this work.

In establishing the optimal tertiary regulation model, the turbine discharges $Q^{\mathrm{turb}}_{\mathrm{H}i}$ at the hydropower plants connect the electromechanical and hydraulic sides.
When pv or load fluctuates, the PI controller in the AGC calculates the incremental power reference for each hydropower plant based on system frequency deviation; see detailed AGC model in literature such as \cite{kundur1994power}.
The AGC participation factors $c_{\mathrm{G}i}$ of the hydropower plants allocate the energy used in the secondary control and thus affect the turbine discharge $Q^{\mathrm{turb}}_{\mathrm{H}i}$ of the plants. Moreover, in islanded operation mode, load can be adjusted by deviating the system frequency reference $\omega^{\mathrm{ref}}$ from zero, and load shedding $P^{\mathrm{shed}}_{\mathrm{D}}$ can be performed in case of insufficient power supply. Consequently, in the tertiary regulation $\bm{c}_{\mathrm{G}}$, $\omega^{\mathrm{ref}}$, and $P^{\mathrm{shed}}_{\mathrm{D}}$ can be adjusted to allocate the water consumption $Q^{\mathrm{turb}}_{\mathrm{H}i}$ among the cascaded hydropower plants to coordinate them.

Compactly, denote the tertiary control variables $\bm{c}_{\mathrm{G}}$, $\omega^{\mathrm{ref}}$ and $P^{\mathrm{shed}}_{\mathrm{D}}$ as a vector,
\begin{align}
 \bm{p} := \left[ \bm{c}_{\mathrm{G}}^{\mathrm{T}}, \omega^{\mathrm{ref}}, P_\mathrm{D}^\mathrm{sh} \right]^{\mathrm{T}}. \label{eq:variables}
\end{align}

Over a time interval $t\in[0,T]$, given time-domain trajectories of the volatile generation of solar or wind energy $\{P_{\mathrm{S}}(t)\}_{t\in[0,T]}$ and load $\{P_{\mathrm{D}}(t)\}_{t\in[0,T]}$, the tertiary control variables $\bm{p}$ determine the electrical power $P_{\mathrm{G}i}$ of each generator, which is denoted by the following function,
\begin{align}
  P_{\mathrm{G}i}(t) &= P_{\mathrm{G}i} \left( P_{\mathrm{S}}(t), P_{\mathrm{D}}(t),\bm{p} \right),\ i\in\Omega^{\mathrm{G}}. \label{eq:orgpg}
\end{align}

With a time horizon compatible with the river dynamics, the mean value of $P_{\mathrm{G}i}(t)$ over the time interval $t\in[0,T]$ is considered, as
\begin{align}
 \hspace{-4pt}\bar{P}_{\mathrm{G}i} \left( \{P_{\mathrm{S}}(t)\}_{t\in[0,T]}, \{P_{\mathrm{D}}(t)\}_{t\in[0,T]},\bm{p} \right) & := \frac{1}{T} \int_{0}^{T} P_{\mathrm{G}i}(t) dt.  \hspace{-4pt} \label{eq:meanpg}
\end{align}

Due to the nonlinearity and complexity of the power system, the analytical expression of (\ref{eq:meanpg}) does not exist. To deal with this problem to establish the system model as Fig. \ref{fig:timescale}(b),  we instead employ an approximate surrogate model to facilitate modeling the optimal model predictive control problem for the tertiary regulation.

First, rewrite the trajectories of the volatile generation and the load as the products of preset normalised trajectories and the corresponding mean values, as
\begin{align}
   P_{\mathrm{S}}(t) & =  \bar{P}_{\mathrm{S}} \times  P^{\mathrm{norm}}_{\mathrm{S}}(t), \
   P_{\mathrm{D}}(t) =  \bar{P}_{\mathrm{D}} \times  P^{\mathrm{norm}}_{\mathrm{D}}(t), \label{eq:loadnorm}
\end{align}
\noindent
with
\begin{align}
   \frac{1}{T}\int_0^T P^{\mathrm{norm}}_{\mathrm{S}}(t) dt = \frac{1}{T} \int_0^T P^{\mathrm{norm}}_{\mathrm{D}}(t) dt = 1;
\end{align}
\noindent
where $\vspace{1pt}\left\{P^{\mathrm{norm}}_{\mathrm{S}}(t)\right\}_{t\in[0,T]}$ and $\left\{P^{\mathrm{norm}}_{\mathrm{D}}(t)\right\}_{t\in[0,T]}$ represent the typical trajectories of the volatile energy and load demand with a unit mean value, respectively.

Substituting (\ref{eq:loadnorm}) into (\ref{eq:orgpg}) and finally into (\ref{eq:meanpg}), we can observe that the mean electrical power of a hydro plant is a function of the mean values of the volatile generation $\bar{P}_{\mathrm{S}}$, the load demand $\bar{P}_{\mathrm{D}}$ over $t\in[0,T]$, and the tertiary control variables $\bm{p}$, denoted as
\begin{align}
  \bar{P}_{\mathrm{G}i} = \bar{P}_{\mathrm{G}i} \left( \bar{P}_{\mathrm{S}}, \bar{P}_{\mathrm{D}} , \bm{p} \right).
\end{align}

Followingly, polynomial approximation is used to construct the surrogate model for $\bar{P}_{\mathrm{G}i}(\cdot)$.
The obtained model has the form as
\begin{align}
  \bar{P}_{\mathrm{G}i} \approx  \bar{P}^{*}_{\mathrm{G}i} \left( \bar{P}_{\mathrm{S}}, \bar{P}_{\mathrm{D}} , \bm{p} \right) :=  \sum\nolimits_{i=1}^{N_\mathrm{b}} f_i \phi_{i} \left( \bar{P}_{\mathrm{S}}, \bar{P}_{\mathrm{D}} , \bm{p} \right), \label{eq:polyapp}
\end{align}
\noindent
where $\phi_{i}(\cdot)$ is the multivariate Legendre polynomial basis function; $f_i$ is the coefficient, obtained by the \emph{collocation method} based on power system dynamic simulation results. See detailed procedure of constructing the polynomial surrogate model in the Appendix A.

To facilitate the formulation of the proposed tertiary regulation as a model predictive controller (MPC) that is solvable as a mixed-integer quadratic programming (MIQP), the approximation order of (\ref{eq:polyapp}) is limited to $2$; see Section \ref{sec:eq} for details. 

Further, generally load shedding is realised by tripping feeders. 
The total amount of load shedding is the sum of the products of the binary variable $\pi_j$ and the capacity of feeders $P_{\mathrm{D}j}^\mathrm{sh}$, as
\begin{align}
   P_\mathrm{D}^\mathrm{sh}(t) = \sum\nolimits_{j \in \Omega^{\mathrm{sh}}} \pi_j P_{\mathrm{D}j}^\mathrm{sh}, \label{eq:shed}
\end{align}
\noindent
and the tertiary control vector is re-denoted as
\begin{align}
 \bm{p} := \left[ \bm{c}_{\mathrm{G}}^{\mathrm{T}}, \omega^{\mathrm{ref}}, \bm{\pi}^{\mathrm{T}} \right]^{\mathrm{T}}. \label{eq:variables}
\end{align}

Substituting (\ref{eq:shed}) into (\ref{eq:polyapp}) and rearranging, we can approximate the mean power of each hydropower generator as a function of $\bar{P}_{\mathrm{S}}$, $\bar{P}_{\mathrm{D}}$, $\bm{c}_{\mathrm{G}}$, $\omega^\mathrm{ref}$, and $\bm{\pi}$, as
\begin{align}
  \bar{P}_{\mathrm{G}i} \approx \bar{P}^*_{\mathrm{G}i}(\bar{P}_\mathrm{S}, \bar{P}_\mathrm{D}, \bm{c}_{\mathrm{G}}, \omega^{\mathrm{ref}}, \bm{\pi}) . \label{eq:pg}
\end{align}

The efficacy of the polynomial surrogate model (\ref{eq:pg}) (or (\ref{eq:polyapp})) is validated numerically. In Section \ref{sec:result}, simulations on PSS/E show that (\ref{eq:pg}) accurately gives the mean values of the hydropower outputs over each $10$-minute interval, which is adequate for establishing a receding-horizon controller.

Finally, we establish the relation between the mean hydropower generation $\bar{P}_{\mathrm{G}i}$ and turbine discharge $Q^{\mathrm{turb}}_{\mathrm{H}i}$ over the time interval $t\in[0,T]$. The electrical power of a hydropower plant is a nonlinear function of water head $ H_{\mathrm{H}i}^\mathrm{head} := H_{\mathrm{H}i}^\mathrm{up} - H_{\mathrm{H}i}^\mathrm{down}$ and $Q^{\mathrm{turb}}_{\mathrm{H}i}$, known as the \emph{production function} \cite{Liu2017interval}, which is expressed as
\begin{align}
  \bar{P}_{\mathrm{G}i} = f_{\mathrm{H}i} \left( H_{\mathrm{H}i}^\mathrm{head}, Q_{\mathrm{H}i}^\mathrm{turb}  \right),\ i\in \Omega^{\mathrm{G}}. \label{eq:product}
\end{align}

To facilitate the formulation of the controller, by linearising (\ref{eq:product}) the turbine discharge is approximated by a linear function as
\begin{align}
  \hat{Q}_{\mathrm{H}i}^\mathrm{turb} = \frac{\bar{P}_{\mathrm{G}i}-\tilde{P}_{\mathrm{G}i}}{\eta_i \rho g  \tilde{H}_{\mathrm{H}i}^\mathrm{head} } - \frac{ \tilde{P}_{\mathrm{G}i}  \hat{H}_{\mathrm{H}i}^\mathrm{head} }{\eta_i \rho g \big( \tilde{H}_{\mathrm{H}i}^\mathrm{head} \big) \hspace{0pt}^2}. \label{eq:linproduct}
\end{align}

Substituting (\ref{eq:pg}) into (\ref{eq:linproduct}), the turbine discharge $\hat{Q}_{\mathrm{H}i}^\mathrm{turb}$ of each hydropower plant during the primary and secondary frequency control process
is finally represented as a polynomial function of of $\hat{H}_{\mathrm{H}i}^\mathrm{up}$, $\hat{H}_{\mathrm{H}i}^\mathrm{down}$, $\bm{c}_{\mathrm{G}}$, $\omega^\mathrm{ref}$, and $\pi_j$. In a word, by adjusting the tertiary control variables (\ref{eq:variables}), the turbine discharge of the hydropower plant over a time interval can be controlled.

Note that the linearisation of (\ref{eq:product}) has also been adopted in many existing works \cite{Hug2011Predictive,qiu2020stochastic}. Because model (\ref{eq:linproduct}) is only used in the tertiary regulation, such linearisation will not affect the primary and secondary control and therefore the frequency stability of the islanded system will not be affected. Further, in online operation the linearisation point is repeatedly updated based on the current operation point, through which the linearisation error can be alleviated in a receding-horizon regulation manner.

\subsection{Hydraulic Side: State-Space Model of River Dynamics}
\label{sec:river}

According to previous researches by the community \cite{Holanda2017Assessment,Hug2011Predictive, Hamann2016Using} as well as our pervious work \cite{qiu2020stochastic}, for modeling the river dynamics in the run-of-the-river hydropower the \emph{shallow water equations} is a suitable model. This section only gives a brief introduction to this model, and details can be found in the above mentioned literature.

The shallow water equations are partial differential equations of the water volume and momentum conservation \cite{Litrico2009}:
\begin{align}
0 & = \frac{\partial Q}{\partial y} + \frac{\partial H}{\partial t}, \label{eq:continuity} \\
0 & = \frac{1}{g} \frac{\partial}{\partial t}\left(\frac{Q}{S}\right) + \frac{1}{2g}\frac{\partial}{\partial y}\left(\frac{Q}{S}\right)^2 + \frac{\partial{H}}{\partial{y}} + I_f - I_0, \label{eq:momentum}
\end{align}
\noindent
where $y$ denotes the position. The equivalent friction slope $I_f$ is empirically modeled by the \emph{Manning-Strickler formula}, which is a function of $Q$ and $H$ (see \cite{Litrico2009}). The parameters such as the river width and slope can be measured in field or estimated via data assimilation \cite{Ding2004Identification}.

\begin{figure}[tb]
	\centering
	\includegraphics[width=2.3in]{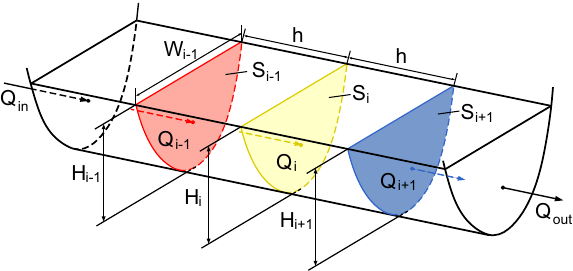}
	\caption{Spatial discretisation framework of the water flow. Detailed formulation can be seen in \cite{qiu2020stochastic}}
	\label{fig:river}
\end{figure}

In normal operation, water stage varies within only a small range compared to the depth. Thus, the river dynamics can also be linearised \cite{Litrico2009}. Then, discretising the water flow into nonoverlapping cells of length $h$ as Fig. \ref{fig:river} and using the finite-difference format, the linear dynamic river model can be obtained \cite{Litrico2009,qiu2020stochastic} compactly as
\begin{align}
{\dot{\bm{x}}_\mathrm{R}(t)} = \bm{A}_\mathrm{R}\bm{x}_\mathrm{R}(t) + \bm{B}_\mathrm{R} \bm{u}_\mathrm{R}\left(t \right) + \bm{C}_\mathrm{R} \bm{z}_\mathrm{R}(t) , \label{eq:river}
\end{align}
\noindent
where $\bm{A}_\mathrm{R}$, $\bm{B}_\mathrm{R}$, and $\bm{C}_\mathrm{R}$ are constant matrices.

The river control variables $\bm{u}_\mathrm{R}\left(t \right)$ including the turbine discharge $\hat{Q}_{\mathrm{H}i}^\mathrm{turb}$ are determined by the operation of the hydropower plants, which are essentially determined by the primary and secondary frequency control and the tertiary regulation. The hydraulic coupling of the cascaded plants is naturally modeled in (\ref{eq:river}). The river operation constraints can also be formed with $\bm{x}_\mathrm{R}$. See detailed deduction of (\ref{eq:river}) in \cite{Hug2011Predictive, Hamann2016Using,qiu2020stochastic}.

In online operation, a state estimator, e.g., the Kalman filter \cite{Glanzmann2005Supervisory}, 
can be employed to provide river state estimation using available measurements such as the turbine discharges and water stages at the plants. This provides full-state feedback for the proposed controller, as shown in Fig. \ref{fig:frame}.

\subsection{The Overall State-Space Model}
\label{sec:overallmodel}

Substituting the second-order polynomial surrogate model (\ref{eq:pg}) of the power system dynamics and the turbine discharge model (\ref{eq:linproduct}) into the dynamic river model (\ref{eq:river}), the overall state-space system model is obtained, denoted compactly as
\begin{align}
  {\dot{\bm{x}}(t)} = \bm{A}\bm{x}(t) & + \bm{B}(\bar{P}_{\mathrm{S}}(t),\bar{P}_{\mathrm{D}}(t)) \bm{u}\left(t \right) \nonumber \\
  & + \bm{C} \bm{z}(t) +  \bm{u}(t)^{\mathrm{T}} \bm{F}\left( \bar{P}_{\mathrm{S}}(t),\bar{P}_{\mathrm{D}}(t) \right) \bm{u}(t) , \label{eq:overall}
\end{align}
\noindent
where $\bm{A}$, $\bm{B}(\cdot)$, $\bm{C}$, and $\bm{F}\left(\cdot\right)$ are coefficient matrices, where $\bm{A}$ and $\bm{C}$ are respectively identical to $\bm{A}_\mathrm{R}$ and $\bm{C}_\mathrm{R}$ in (\ref{eq:river}); $\bm{x}(t)$ and $\bm{z}(t)$ are the same as $\bm{x}_\mathrm{R}(t)$ and $\bm{z}_\mathrm{R}(t)$; 
\begin{align}
 \bm{u}(t):=\left[\bm{c}_{\mathrm{G}}(t)^\mathrm{T}, \omega^{\mathrm{ref}}(t), \bm{\pi}(t)^\mathrm{T}, \bm{Q}^\mathrm{sp}_\mathrm{H}(t)^\mathrm{T} \right]^\mathrm{T}
\end{align}
\noindent
is the vector of decision variables of the tertiary regulation; the quadratic term $\bm{u}(t)^{\mathrm{T}} \bm{F}\left( P_{\mathrm{D}}(t) \right) \bm{u}(t)$ represents the quadratic terms of control variables in the second-order polynomial surrogate model (\ref{eq:pg}); and $\bm{B}(\bar{P}_{\mathrm{S}}(t),\bar{P}_{\mathrm{D}}(t))$ and $\bm{F}\left( \bar{P}_{\mathrm{D}}(t) \right)$ indicate these matrices are functions of $\bar{P}_{\mathrm{S}}(t)$ and $\bar{P}_{\mathrm{D}}(t)$.

The framework of this model is illustrated in Fig. \ref{fig:timescale}(b).

\subsection{Discussion on the Dimensionality of Modeling}
\label{sec:dimension}

In the proposed model (\ref{eq:overall}), because the frequency dynamics are modeled by algebraic functions (\ref{eq:pg}), the state variables are related only to the river dynamics. In the case study in Section \ref{sec:case}, the size of the dynamic river model (\ref{eq:river}) is $200$. Considering the MPC step length $T = 600$ s and horizon $N = 12$ or $7,200$ s, the size of the discrete system model (\ref{eq:eqstate}) is $200 \times 12 = 2400$, which is appropriate for online rolling optimisation.

Otherwise, if additionally considering the detailed power system dynamics with total $n$ ($n$ varies from tens to thousands depending on system complexity) state and algebraic variables, with a step length of $0.02$ s that is compatible with the electromechanical dynamics, the overall size of the state reaches $(n + 200) \times 50 \times 7200$. Even if the power system and river dynamics are discretised in different time resolutions, the size of the states reaches $n \times 50 \times 7200 + 12 \times 200$, much larger than the proposed model. Considering using interior point method (IPM) with time complexity $O(n^4)$ to solve the MPC, the efficiency of the proposed model is significantly better.

\section{Mathematical Formulation of the Coordinated Tertiary Regulation}
\label{sec:controller}

As shown in Fig. \ref{fig:frame}, for the tertiary regulation an MPC is employed to repeatedly optimise the tertiary control variables including the AGC participation factors to coordinate the cascaded plants and to give commands of the frequency reference and load shedding in a receding-horizon manner.

There are two main reasons for this work to adopt the MPC as the tertiary controller. First, the cascaded run-of-the-river hydropower system is a dynamic system modeled by the state space equation (\ref{eq:overall}) with various operational inequality constraints; MPC is very suitable for the optimal control of such a constrained dynamic system. Second, due to the stochastic nature of pv, the operation trajectory of the system may deviate from the predicted one, which can be well alleviated by the well-known receding horizon scheme of the MPC.

Given the prediction horizon $N$ and step length $T$, the objective and constraints of the MPC are formulated, as explained below.

\subsection{Objective Function}
\label{sec:objective}

The overall control objective involves coordinating the cascaded plants and reducing load shedding. Meanwhile, the river state, such as the water stage, should not deviate too far from the nominal. The detailed objective includes the following components.

\vskip 6pt
{\bf Deviation in the Frequency Reference:}
According to (\ref{eq:pg}), the load power can be adjusted by the frequency reference, which assists power balance in the islanded system. However, the frequency bias should not deviate from zero if unnecessary.
Therefore, the following quadratic function is included in the objective:
\begin{align}
  J_1 = \sum\nolimits_{k=0}^{N-1} \omega^\mathrm{ref}(kT)^2. \label{eq:objfreq}
\end{align}

\vskip 6pt
{\bf Load Shedding:}
On the premise of power balance and frequency stability, load shedding should be minimised to avoid loss of load, which is expressed as:
\begin{align}
  J_2 = \sum\nolimits_{k=0}^{N-1} P_\mathrm{D}^\mathrm{sh}(kT). \label{eq:objshed}
\end{align}

\vskip 6pt
{\bf River Stage Deviation:}
During the regulation process, the water stage and discharge along the river and channels of the plants should not deviate far from the nominal. This can be achieved by minimising the following quadratic function:
\begin{align}
  J_3 = \sum\nolimits_{k=1}^{N} \bm{x}_\mathrm{R}(kT)^\mathrm{T}\bm{x}_\mathrm{R}(kT). \label{eq:objriver}
\end{align}

\vskip 6pt
{\bf Water Spillage:}
When the upstream inflow exceeds the power demand and the upper limit of the water stage is encountered, plant water spillage is needed to ensure operational security. Spillage occurs only when needed, represented as minimising
\begin{align}
  J_4 = \sum\nolimits_{k=0}^{N-1} \sum\nolimits_{j \in \Omega^{\mathrm{G}}} Q_{\mathrm{G}i}^\mathrm{sp}(kT).
\end{align}

\vskip 6pt
{\bf Quadratic Terms of the AGC Participation Factors:}
Finally, to avoid oscillation in the AGC participation factors, the following quadratic term is included in the objective:
\begin{align}
  J_5 = \sum\nolimits_{k=0}^{N-1}  \bm{c}_{\mathrm{G}}(kT)^\mathrm{T} \bm{c}_{\mathrm{G}}(kT).
\end{align}

The overall control objective is defined as a weighted sum of the above terms with positive weight parameters:
\begin{align}
  J = \lambda_1 J_1 + \lambda_2 J_2 + \lambda_3 J_3 + \lambda_4 J_4 + \lambda_5 J_5. \label{eq:obj}
\end{align}
\noindent

The weights in the objective function (\ref{eq:obj}) are chosen based on a tradeoff between different components. For example, when power supply and reliability are more preferred, larger $J_1$ and $J_2$ can be chosen; and when river operational constraints are more critical, $J_3$ should be increased. In this paper, we set $\lambda_1=10$, $\lambda_2=10$, $\lambda_3=1$, $\lambda_4=1,000$, and $\lambda_5=10$. In Section \ref{sec:weights}, the impact of these weights on control performance is numerically demonstrated.

\subsection{Equality Constraints}
\label{sec:eq}

The equality constraints in the proposed controller include two components, as listed below.

\vskip 6pt
{\bf System Dynamics Model:}
The state-space model (\ref{eq:overall}) is temporally discretised into the equality constraints as
\begin{align}
  \bm{x}\big( (k+1) &  T\big) = \bar{\bm{A}} \bm{x}(kT)  + \bar{\bm{B}}\big( \bar{P}_{\mathrm{S}}(kT),\bar{P}_{\mathrm{D}}(kT) \big) \bm{u}\left(kT \right) \nonumber \\
  +   \bar{\bm{C}} & \bm{z}(kT) + \bm{u}(kT)^{\mathrm{T}} \bar{\bm{F}} \big(  \bar{P}_{\mathrm{S}}(kT), \bar{P}_{\mathrm{D}}(kT) \big) \bm{u}(kT) \label{eq:eqstate}
\end{align}
\noindent
for $k=0,\ldots,N-1$, where $\bar{\bm{A}}$, $\bar{\bm{B}}$, $\bar{\bm{C}}$ and $\bar{\bm{F}}$ are coefficient matrices;
the solar/wind power $\bar{P}_{\mathrm{S}}(kT)$, load power $\bar{P}_{\mathrm{D}}(kT)$, and river inflow $\bm{z}(kT)$ in the future are given as forecasts by the dispatching system, and are constant parameters in the optimisation problem.

When using second-order polynomials to approximate the mean power output of the hydropower plants, as introduced in Section \ref{sec:modelelectric} and exemplified in Table \ref{tab:poly}, bilinear terms in terms of the tertiary control variables, i.e., $c_{\mathrm{G}i} \omega^\mathrm{ref}$ and $\pi_j c_{\mathrm{G}i}$  are included as the last term of (\ref{eq:eqstate}). Directly incorporating them into an optimisation may cause non-convexity, making the optimal tertiary regulation problem difficult to solve.

 To address the bilinear terms to make the optimisation problem convex, $\omega^\mathrm{ref}$ is discretised as
\begin{align}
   \omega^\mathrm{ref} =  \underline{\omega}^\mathrm{ref} + \sum\nolimits_{j=0}^{K-1} 2^j \Delta \omega  \beta_j,
\end{align}
\noindent
where $\Delta\omega = (\overline{\omega}^\mathrm{ref} - \underline{\omega}^\mathrm{ref})/{2^M}$ is the discretisation step length; $K>0$ is an integer; $\beta_j$ is a binary variable. Thus, the bilinear term $c_{\mathrm{G}i} \omega^\mathrm{ref}$ is replaced with the linear combination of $c_{\mathrm{G}i}\beta_j$. 

Then, the big-$M$ method is used to address $c_{\mathrm{G}i} \beta_j$ as
\begin{align}
  c_{\mathrm{G}i} \omega^\mathrm{ref}  =   c_{\mathrm{G}i}\underline{\omega}^\mathrm{ref}  & + \sum\nolimits_{j=0}^{K-1} \Delta\omega 2^j \delta^{(\omega)}_{i,j},  \label{eq:comega} \\
  \Delta\omega - M (1-\beta_j)  \le & \delta^{(\omega)}_{i,j} \le \Delta\omega + M (1-\beta_j), \label{eq:comega1}\\
  - M \beta_j  \le &\delta^{(\omega)}_{i,j} \le M \beta_j. \label{eq:comega2}
\end{align}

Similarly, the bilinear term $c_{\mathrm{G}i} \pi_j$ is addressed as
\begin{align}
  c_{\mathrm{G}i} \pi_j  & =  \delta^{(c_{\mathrm{G}i})}_{i,j}, \label{eq:cpi}\\
   c_{\mathrm{G}i} - M (1-\pi_j)  \le & \delta^{(c_{\mathrm{G}i})}_{i,j}  \le  c_{\mathrm{G}i} + M (1-\pi_j), \label{eq:cpi1}\\
  - M \pi_j  \le & \delta^{(c_{\mathrm{G}i})}_{i,j} \le  M \pi_j. \label{eq:cpi2}
\end{align}

Substituting (\ref{eq:comega}) and (\ref{eq:cpi}) into (\ref{eq:eqstate}), the bilinear terms are replaced with linear combinations of the intermediate variables $\delta^{(\omega)}_{i,j}$ and $\delta^{(c_{\mathrm{G}i})}_{i,j}$. Thus far, a convex programming model can be formulated with linear constraints (\ref{eq:comega1})--(\ref{eq:comega2}) and (\ref{eq:cpi1})--(\ref{eq:cpi2}).

\vskip 6pt
{\bf Sum of the AGC Participation Factors:}
Obviously, the AGC participation factors sum up to $1$, as 
\begin{align}
  1 = \sum\nolimits_{i \in \Omega^{\mathrm{G}}}  c_{\mathrm{G}i}(kT), \ k=0,\ldots, N-1.
\end{align}

\subsection{Inequality Constraints}
\label{sec:ieq}

The inequality constraints for $k=0,\ldots, N-1$ include the following components.

\vskip 6pt
{\bf AGC Participation Factors:} The AGC participation factors should be within the interval $[0,1]$ as
\begin{align}
   0 \le c_{\mathrm{G}i}(kT) \le 1,\ i \in \Omega^{\mathrm{G}}.
\end{align}

\vskip 6pt
{\bf Hydropower Generation Limits:} The power output of each hydropower plant should not exceed the limits, as follows:
\begin{align}
   \underline{P}_{\mathrm{G}i} \le  \bar{P}_{\mathrm{G}i}(kT;c_{\mathrm{G}i}, \omega^{\mathrm{ref}}, \bm{\pi}) \le \overline{P}_{\mathrm{G}i}, \ i \in \Omega^{\mathrm{G}}.
\end{align}
The bilinear terms in $ \hat{P}_{\mathrm{G}i}(\cdot)$ are replaced with (\ref{eq:comega}) and (\ref{eq:cpi}), and the resulting constraints are linear.

\vskip 6pt
{\bf Electrical Network Constraints:} Network constraints such as the branch limits can be considered as dc power flow:
\begin{align}
   \underline{\bm{P}}_\mathrm{B} \le \bm{A}_\mathrm{Net} \bm{B}_\mathrm{Net}^{-1} \bm{P}_{\mathrm{Inj}} \le \overline{\bm{P}}_\mathrm{B}.
\end{align}

\vskip 6pt
{\bf Water State Limits:} The operation regulations require that the water stage stay within an allowed interval. This is considered at monitoring points along the river and channels:
\begin{align}
\underline{H}_{i} \le H_{i}(t) \le \overline{H}_{i},\ i \in \Omega^\mathrm{m}. \label{eq:stagelimit}
\end{align}

If additional constraints such as the 
switching counting limit of load shedding need to be considered, they can also be easily included in the MPC. To save space, this is not discussed here.

\subsection{MPC Formulation}
\label{sec:overallctrl}

Summarising all the above, the overall MPC problem in the proposed coordinated tertiary regulation is established, denoted as
\begin{align}
  \min_{\bm{U}} \ (\ref{eq:obj}),  \  \text{subject to}  \hspace{6pt} (\ref{eq:eqstate})-&(\ref{eq:stagelimit}), \label{eq:all}
\end{align}
\noindent
where $\bm{U} = [ \bm{u}(0)^\mathrm{T},\ldots,\bm{u}\big((k-1)T\big)\hspace{0pt}^\mathrm{T} ]\hspace{0pt}^\mathrm{T}$ is the sequence of the control variables.
The MPC problem (\ref{eq:all}) is a typical mixed-integer quadratic programming (MIQP), which can be solved using commercial solvers such as \emph{IBM ILOG Cplex}.

Once the MPC problem is solved, the AGC participation factors, frequency reference, and load shedding are updated based on the first entry of $\bm{U}$. 
When shifted to the next step, the MPC is solved again based on the updated state estimation and forecasts in a receding-horizon manner.
The power references of the plants are updated in real time by the PI controller and the repeatedly updated participation factors of the AGC, shown in Fig. \ref{fig:frame}.

Note that although modeling errors are introduced in approximating the power system dynamics in Section \ref{sec:modelelectric} and linearising the river dynamics in Section \ref{sec:river}, they can be alleviated by the receding-horizon scheme.

\begin{figure}[tb]
	\centering
	\includegraphics[width=2.8in]{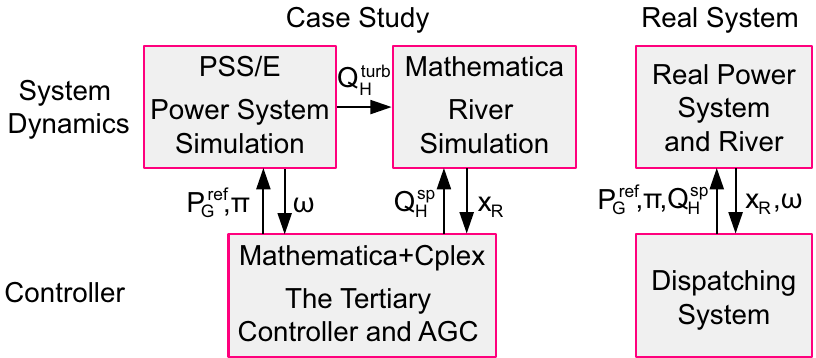}
	\caption{The simulation framework compared with the actual system}
	\label{fig:simulation}
\end{figure}

{\color{black}
\subsection{Online Smoothing Control of PV to Reduce Volatility}
\label{sec:smoothing}

As noted in Introduction, with a high pv penetration, the stochastic volatility may exceed the adjustability of the primary and secondary frequency control of the hydropower units and therefore poses an adverse impact on stability. Hence, before applying the above proposed tertiary control, a smoothing control of pv is needed to reduce the impact of stochastic variation on frequency stability.

Supposing if the pv power smoothly follows its prediction trajectory, its impact on system stability can be well counterbalanced. However, we cannot simply track the prediction because the available pv power may drop below it. Hence, a controlled curtailment that considers the real-time volatility of pv is adopted in this work.

Specifically, the controlled curtailment of pv at minute $k$ is as,
\begin{align}
  P_{\mathrm{S}}(k) & = \min \big\{ P_{\mathrm{S}}^{\mathrm{pred}}(k) - P_{\mathrm{S}}^{\mathrm{cur}}(k)(\ge0), P_{\mathrm{S}}^{\mathrm{avail}}(k) \big\}, \label{eq:cur} \\
  P_{\mathrm{S}}^{\mathrm{cur}}(k) & = \alpha \kappa_\gamma  v(k) + (1-\alpha)P_{\mathrm{S}}^{\mathrm{cur}}(k-1) \label{eq:smooth} \\
  v(k) &=  \sqrt{ \frac{1}{\tau} \sum_{t=k-\tau}^{k-1} \big[ P_{\mathrm{S}}^{\mathrm{pred}}(t) - P_{\mathrm{S}}^{\mathrm{avail}}(t)\big]^2 }  \label{eq:var}
\end{align}
\noindent
where $P_{\mathrm{S}}^{\mathrm{pred}}(k)$ is the prediction and $P_{\mathrm{S}}^{\mathrm{cur}}(k)$ is the curtailment at minute $k$; $P_{\mathrm{S}}(k)$ is the smoothed output; $P_{\mathrm{S}}^{\mathrm{avail}}(k)$ is the available pv power determined by irradiance, usually the maximum power point (MPP); (\ref{eq:var}) calculates the volatility of pv for the past $\tau$ minutes; (\ref{eq:smooth}) is an low-pass filter with parameter $\alpha$ to make the curtailment varies smoothly; $\kappa_\gamma$ is the level of curtailment.

When irradiance varies randomly to cause the available pv power varies significantly, $v(k)$ becomes large and the curtailment of pv is increased to counter the volatility. In contrast, when pv volatility is mild, $v(k)$ becomes small and the pv output largely follows the prediction. Thus, the output of pv becomes relatively smooth under different circumstance. Examples are illustrated in Section \ref{sec:smoothingresult}.

}

\section{Case Study}
\label{sec:case}

\subsection{Simulation Platform}

To verify the proposed regulation approach, a detailed simulation platform is established jointly based on \emph{PTI PSS/E 34} and \emph{Wolfram Mathematica 11.3}, as shown in Fig. \ref{fig:simulation}. The electrical side is based on PSS/E, including electrical network, GENROU generator, SCRX exciter, PSS2A stabiliser, and HYGOVM governor models which include detailed penstock, turbine, and governor dynamics. The river model is based on \emph{Mathematica}, including the shallow-water-equation-based river model (\ref{eq:continuity})--(\ref{eq:momentum}), solved by the finite difference method. The electrical and hydraulic sides communicate with each other via the \emph{PSSPY} interface. The proposed controller (\ref{eq:all}) is modeled on \emph{Mathematica} and solved by \emph{IBM ILOG Cplex} via the \emph{NETLink} interface. The step length of dynamic river simulation is $10$ s, and the control period of AGC is $4$ s.

\subsection{Case Setting}

\begin{figure}[tb]
   { \centering
    \subfloat[]{\includegraphics[width=2.3in]{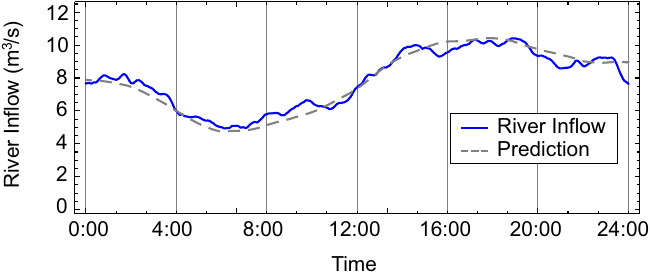}} \\
    \subfloat[]{\includegraphics[width=2.3in]{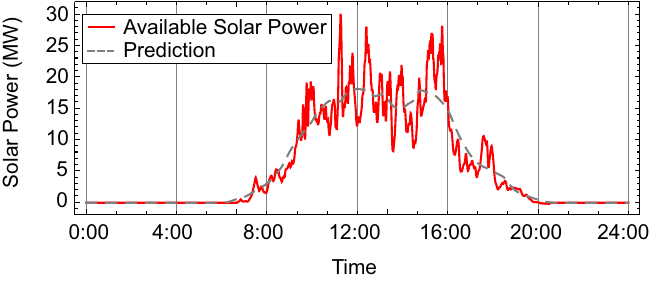}} \\
    \subfloat[]{\includegraphics[width=2.3in]{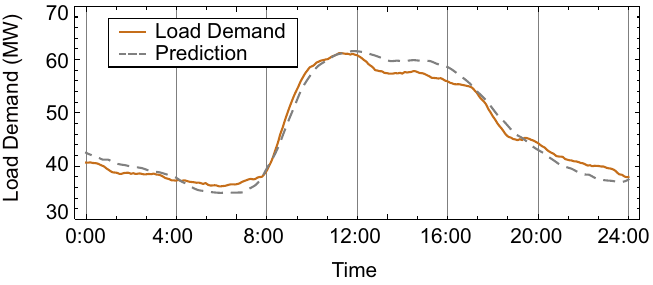}} \\}
  \caption{Simulation scenario used in case study in Sections \ref{sec:polyresult} to \ref{sec:result}}
    \figfooter{a}{River upstream inflow and the forecast}
    \figfooter{b}{{\color{black}Available} solar power and the forecast} 
    \figfooter{c}{Load and the forecast}
  \label{fig:scen}
\end{figure}

The real-life system located in Xiaojin County, Sichuan Province, China, shown in Fig. \ref{fig:xiaojin}, is used in this case study.
{\color{black}
Modeling parameters for simulation are specified in Appendix B.}
Specifically, data for the ten river sections are listed in Table \ref{tab:river}, and brief data for the three cascaded hydropower plants (MP, YJW, and MGQ, from upstream to downstream) that participate in AGC are given in Table \ref{tab:plant}. The hydro turbines at the three cascaded plants are of Francis type, and the detailed parameters of the HYGOVM models for the turbine-governors used in simulation are given in Table \ref{tab:turbine}. Other four smaller hydropower plants (HJQ, HK, SGQ, and MW, at bus 34, 50, 44, and 42, respectively; see Fig. \ref{fig:xiaojin}) do not participate in AGC. Their power references are fixed at $2.60$, $2.70$, $4.00$, and $4.20$ MW, respectively.
{\color{black} The deadband of the primary control is set to be $\pm0.05$ Hz around the setting point in such a small system. The electrical network parameters are given in Table \ref{tab:netpara}.}
The parameters specified in Tables \ref{tab:river}, \ref{tab:plant}, \ref{tab:turbine}, and \ref{tab:netpara} are based on on-site collected data, which are regarded as constants in the proposed control scheme.

The ramping ability of the hydropower plants is limited by the governor time constant, gate opening/closing rate, and surge chamber and penstock dynamics specified in Table \ref{tab:turbine}.
Five feeders of power $1$, $1$, $2$, $2$, and $4$ MW at buses 36 and 37 serve as sheddable loads, and each one is permitted to switch once per hour. 

We choose to test the proposed regulation method in the dry season, as in the wet season, there is always abundant water to generate enough electrical power. In contrast, in the dry season, water resources are limited, and the total hydropower generation cannot satisfy the load demands without solar power generation and load shedding. In this situation, the adjustability of the hydropower should be fully exploited, where the proposed coordinated regulation shows its value.

Specifically, for this case study, the upstream water inflow and its forecast are as shown in Fig. \ref{fig:scen}(a), being much smaller than the rated turbine discharge of $30$ to $40$ $\text{m}^3/\text{s}$ of the whole plants each with $3$ turbines. Thus, of each plant only one turbine is in operation. The pv and load power and their forecasts are respectively shown in Figs. \ref{fig:scen}(b) and \ref{fig:scen}(c). Specifically, the pv data in Fig. \ref{fig:scen}(b) is the real-time pv power collected in the Xiaojin system on May 30, 2018.

By examining historical data, we find the maximal volatility of pv is $11.59$ MW/min, and in $99.9\%$ of the operation time pv volatility is within $\pm6.99$ MW/min.
{\color{black} To alleviate such random volatility, the smoothing control introduced in Section \ref{sec:smoothing} is adopted with result is demonstrated in Section \ref{sec:smoothingresult}.
}

In the MPC of the tertiary regulation, the prediction step length is set as $T=600$ s. Since the water wave travels through the cascaded plants for one more hour, we set the prediction horizon of the controller to be two-times longer, i.e., $2$ hour or $N=12$. The objective function is set as (\ref{eq:obj}) with $\lambda_1=10$, $\lambda_2=10$, $\lambda_3=1$, $\lambda_4=1,000$, and $\lambda_5=10$. The river operational constraints include water stage limits at the monitoring points $800$ m upstream of the dams and plants. On natural river reaches and channels, the limits are $\pm 0.2$ m and $\pm 0.5$ m around the nominal point, respectively. The frequency reference limit is set as $\pm 0.1$ Hz.

{\color{black}

\subsection{Result of PV Smoothing Control}
\label{sec:smoothingresult}

\begin{figure}[tb]
   { \centering
    \subfloat[]{\includegraphics[width=2.3in]{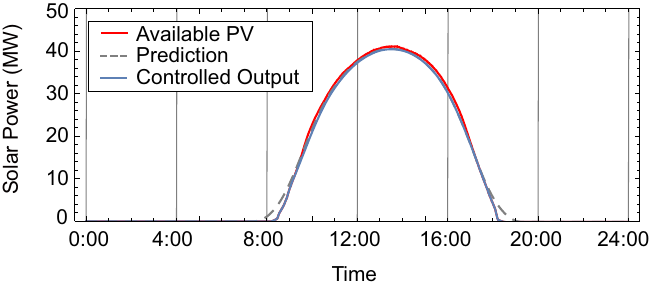}} \\
    \subfloat[]{\includegraphics[width=2.3in]{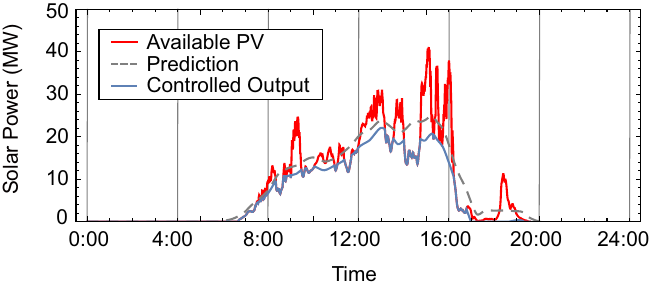}} \\
    \subfloat[]{\includegraphics[width=2.3in]{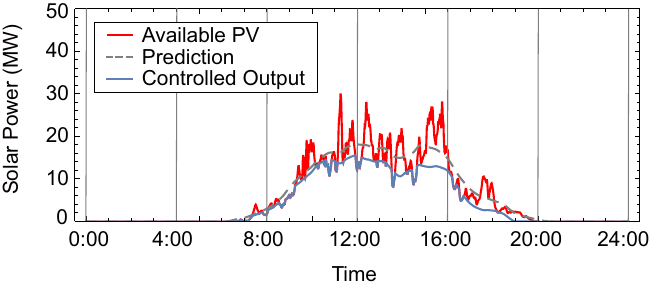}} \\}
  \caption{\color{black}Available pv power and simulation result of smoothing control output on typical days}\color{black}
    \figfooter{a}{Jan 10, 2018}
    \figfooter{b}{May 5, 2018}
    \figfooter{c}{May 30, 2018, which used for simulation in Sections \ref{sec:benchmark}--\ref{sec:discussion}  }
  \label{fig:smooth}
\end{figure}

\begin{figure}[tb]
	\centering
	\includegraphics[width=3.3in]{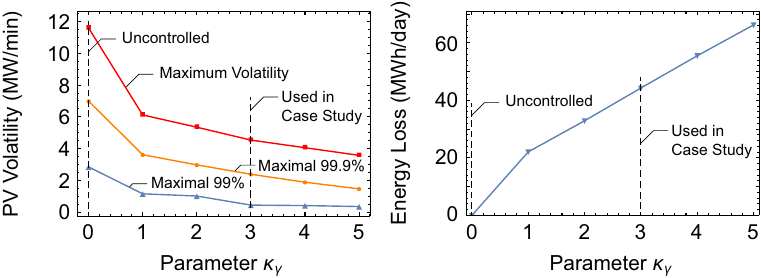}
	\caption{\color{black}Volatility of pv generation under smoothing control and the average daily loss of energy with different parameter $\kappa_\gamma$}
	\label{fig:smoothpara}
\end{figure}

Simulation result of the smoothing control introduced in Section \ref{sec:smoothing} with $\tau=15$ and $\kappa_\gamma=3$ is illustrated on different days of available pv power collected in the Xiaojin system shown in Fig. \ref{fig:smooth}. As can be seen, when the volatility is mild as shown in Fig. \ref{fig:smooth}(a), the curtailment is very small, and the pv output largely follows the prediction. In contrast, when intensive volatility occurs as Figs. \ref{fig:smooth}(b) and \ref{fig:smooth}(c), the curtailment of pv is controlled according to the intensity of volatility, and the output follows a similar trend to the prediction but in a lower profile, which significantly reduces the variations of pv and therefore alleviates the challenge of frequency stability.

Further, the smoothing effect and energy loss due to curtailment control are investigated based on the pv data for the whole year of 2018. The maximum, and $99.9\%$ and $99\%$ largest minutely volatility of pv, and the average daily energy loss with different $\kappa_\gamma$ are shown in Fig. \ref{fig:smoothpara}. Here, $\kappa_\gamma=0$ means no smoothing control. Given larger $\kappa_\gamma$, the pv output is smoother but energy loss becomes larger. In practice, $\kappa_\gamma$ can be flexibly selected based on the robustness of the power system. In the case study we choose $\kappa_\gamma=3$.

}

\subsection{Polynomial Surrogate Model for the Power System Dynamics}
\label{sec:polyresult}

Using the polynomial approximation method presented in the Appendix A with the error tolerance $\epsilon=0.01$ in the Smolyak approximation algorithm, polynomial surrogate models to characterise the mean power of the three cascade plants over a time interval of $10$ minutes are obtained. The normalised trajectories of solar power and load demand, i.e., $P^{\mathrm{norm}}_{\mathrm{S}}(t)$ and $P^{\mathrm{norm}}_{\mathrm{D}}(t)$ in (\ref{eq:loadnorm}) are based on the ones shown Figs. \ref{fig:scen}(b) and \ref{fig:scen}(c) from 12:00 to 12:10. The obtained models are given in Table \ref{tab:poly}.

As can be seen, the quadratic terms with respect to the control variables in the polynomial functions include $c_{\mathrm{G}i}\omega^{\mathrm{ref}}$ and $c_{\mathrm{G}i}P_{\mathrm{D}}^{\mathrm{sh}}$ (or $c_{\mathrm{G}i}\pi_{j}$ after substituting (\ref{eq:shed}) in). By using the transformation (\ref{eq:comega})--(\ref{eq:cpi}) introduced in Section \ref{sec:eq}, these quadratic terms can be replaced by linear terms, making the optimal tertiary regulation model (\ref{eq:all}) an MIQP.

\begin{figure}[tb]
	\centering
	\includegraphics[width=2.3in]{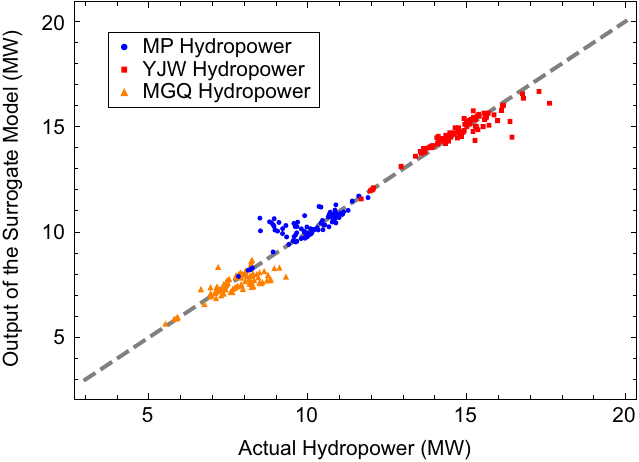}
	\caption{\color{black}Output of the polynomial surrogate models shown in Table \ref{tab:poly} versus the actual generation of the three cascaded plants, where each mark represents the mean generation over a 10-min interval}
	\label{fig:poly}
\end{figure}

The polynomial surrogate models are verified numerically. After the simulation later performed in Section \ref{sec:result}, the actual power outputs of the three cascaded plants and the values obtained by the polynomial surrogate models are compared in Fig. \ref{fig:poly}. We can see that the marks lie close to the line $y=x$, indicating that the employed surrogate models have good accuracy. Further, the decent controller performance shown in Section \ref{sec:result} again verifies these models.

Moreover, we calculate the surrogate models using different normalised trajectories of $P^{\mathrm{norm}}_{\mathrm{S}}(t)$ and $P^{\mathrm{norm}}_{\mathrm{D}}(t)$, and the accuracy of the obtained models are almost the same and good as well. 

\begin{table}[tb]\scriptsize
	\processtable{Polynomial surrogate models for the mean power of the cascaded plants, which are used in the proposed tertiary regulation in the case study\label{tab:poly}}
	{\begin{tabular}{@{\extracolsep{\fill}}ll}
		\toprule
		Hydropower Plant                                      & \phantom{$+$}Polynomial Surrogate Model     \\
		\midrule
		\#1 MP \hspace{1.5em} $\bar{\text{P}}_\text{G1}$=
                     & \phantom{$+$}8.441$+$0.319 $\bar{\text{P}}_\text{D}$$-$0.311$\bar{\text{P}}_\text{S}$$+$0.000149$\bar{\text{P}}_\text{D}$$\bar{\text{P}}_\text{S}$ \\
                     & $-$(0.0174$-$0.00671$\bar{\text{P}}_\text{D}$$+$0.00719$\bar{\text{P}}_\text{S}$)c$_{\text{G1}}$\\
                     & $+$(0.00573$-$0.00191$\bar{\text{P}}_\text{D}$$+$0.00196$\bar{\text{P}}_\text{S}$)c$_{\text{G2}}$\\
                     & $+$(0.0130$-$0.00437$\bar{\text{P}}_\text{D}$$+$0.00456$\bar{\text{P}}_\text{S}$)c$_{\text{G3}}$  \\
                     & $+$(5.551$+$0.00571$\bar{\text{P}}_\text{D}$$-$0.00194$\bar{\text{P}}_\text{S}$)\textomega$^{\text{ref}}$ \\
                     & $-$(0.309$+$0.000533$\bar{\text{P}}_\text{D}$$+$0.000143$\bar{\text{P}}_\text{S}$)P$_\text{D}^{\text{sh}}$ \\
                     & $+$0.127c$_{\text{G1}}$\textomega$^{\text{ref}}$$-$0.0405c$_{\text{G2}}$\textomega$^{\text{ref}}$$-$0.0731c$_{\text{G3}}$\textomega$^{\text{ref}}$ \\
                     & $-$0.00543c$_{\text{G1}}$P$_\text{D}^{\text{sh}}$$+$0.00202c$_{\text{G2}}$P$_\text{D}^{\text{sh}}$$+$0.00368c$_{\text{G3}}$P$_\text{D}^{\text{sh}}$ \\
        \midrule
        \#2 YJW \hspace{1.5em} $\bar{\text{P}}_\text{G2}$=
                     & \phantom{$+$}12.3$+$0.428$\bar{\text{P}}_\text{D}$$-$0.416$\bar{\text{P}}_\text{S}$$+$0.000198 $\bar{\text{P}}_\text{D}$$\bar{\text{P}}_\text{S}$ \\
                     & $+$(0.00841$-$0.00372$\bar{\text{P}}_\text{D}$$+$0.00376$\bar{\text{P}}_\text{S}$)c$_{\text{G1}}$ \\
                     & $-$(0.00830$-$0.00355$\bar{\text{P}}_\text{D}$$+$0.00374$\bar{\text{P}}_\text{S}$)c$_{\text{G2}}$ \\
                     & $+$(0.0138$-$0.00562$\bar{\text{P}}_\text{D}$$+$0.00600$\bar{\text{P}}_\text{S}$)c$_{\text{G3}}$ \\
                     & $+$(7.508$+$0.00747$\bar{\text{P}}_\text{D}$$-$0.00292$\bar{\text{P}}_\text{S}$)\textomega$^{\text{ref}}$ \\
                     & $-$(0.414$+$0.000712$\bar{\text{P}}_\text{D}$$+$0.000182$\bar{\text{P}}_\text{S}$)P$_\text{D}^{\text{sh}}$\\
                     & $-$0.104c$_{\text{G1}}$\textomega$^{\text{ref}}$$+$0.0731 c$_{\text{G2}}$\textomega$^{\text{ref}}$$-$0.0909c$_{\text{G3}}$\textomega$^{\text{ref}}$ \\
                     & $+$0.00446c$_{\text{G1}}$P$_\text{D}^{\text{sh}}$$-$0.00343c$_{\text{G2}}$P$_\text{D}^{\text{sh}}$$+$0.00509 c$_{\text{G3}}$P$_\text{D}^{\text{sh}}$ \\
        \midrule
        \#3 MGQ \hspace{1.5em} $\bar{\text{P}}_\text{G3}$=
                     & \phantom{$+$}6.115$+$0.253$\bar{\text{P}}_\text{D}$$-$0.247$\bar{\text{P}}_\text{S}$$+$0.000118$\bar{\text{P}}_\text{D}$$\bar{\text{P}}_\text{S}$ \\
                     & $+$(0.00476$-$0.00252$\bar{\text{P}}_\text{D}$$+$0.00258$\bar{\text{P}}_\text{S}$)c$_{\text{G1}}$ \\
                     & $+$(0.00386$-$0.00172$\bar{\text{P}}_\text{D}$$+$0.00182$\bar{\text{P}}_\text{S}$)c$_{\text{G2}}$ \\
                     & $-$(0.0249$-$0.00980$\bar{\text{P}}_\text{D}$$-$0.0103$\bar{\text{P}}_\text{S}$)c$_{\text{G3}}$ \\
                     & $+$(4.401$+$0.00454$\bar{\text{P}}_\text{D}$$-$0.00180$\bar{\text{P}}_\text{S}$)\textomega$^{\text{ref}}$ \\
                     & $-$(0.247$+$0.000423$\bar{\text{P}}_\text{D}$$+$0.000103$\bar{\text{P}}_\text{S}$)P$_\text{D}^{\text{sh}}$ \\
                     & $-$0.0561c$_{\text{G1}}$\textomega$^{\text{ref}}$$-$0.0347c$_{\text{G2}}$\textomega$^{\text{ref}}$$+$0.195c$_{\text{G3}}$\textomega$^{\text{ref}}$ \\
                     & $+$0.00281c$_{\text{G1}}$P$_\text{D}^{\text{sh}}$$+$0.0017c$_{\text{G2}}$P$_\text{D}^{\text{sh}}$$-$0.0100c$_{\text{G3}}$P$_\text{D}^{\text{sh}}$ \\
		\botrule
	\end{tabular}}{}
\end{table}

\subsection{Benchmarking Regulation Methods and Simulation}
\label{sec:benchmark}

Four benchmarking regulation methods (denoted as BMs hereafter) are used for comparison:

{\bf BM1:}
PI-based AGC with fixed participation factors that are proportional to the capacities of the cascaded plants, $\omega^{\mathrm{ref}}=0$, and no load shedding.

\begin{figure}[tb]
	\centering
	\includegraphics[width=2.7in]{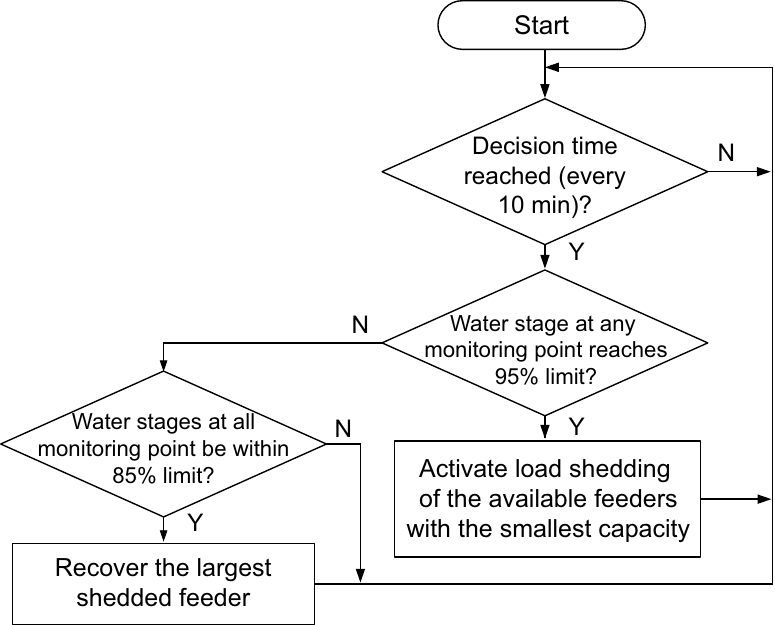}
	\caption{Logic diagram of load shedding control in BM2 and BM3}
	\label{fig:logic}
\end{figure}

{\bf BM2:}
PI-based AGC with fixed participation factors that are proportional to the capacities of the cascaded plants, and $\omega^{\mathrm{ref}}=0$. Load shedding is based on a mixed-logic controller, as shown in Fig. \ref{fig:logic}. The decision period of load shedding is $10$ min.

{\bf BM3:}
Mostly the same as BM2, expect that the frequency reference is set as $\omega^{\mathrm{ref}}=-0.1$ Hz to minimise load energy consumption and to minimise load shedding as the result.

{\bf BM4:}
A receding-horizon generation scheduling with $N=12$ and $T=600$ s is used to optimise the nominal power references of the hydropower plants, the frequency reference, and the load shedding every $10$ minutes. The traditional PI-based AGC with fixed participation factors calculates the incremental power references to stabilise frequency in real time within each period of $10$ minutes. The scheduling model is similar to (\ref{eq:all}) but without consideration of the adjustment of the AGC participation factors.

The benchmarking methods BM1 to BM3 are based on common idea of power system regulation. Specially note that optimisation-based BM4 is also first proposed in this paper. As briefly mentioned in the Introduction, to the best of the authors' knowledge, there is no publication on the tertiary regulation of islanded power systems with cascaded run-of-the-river hydropower. Due to space limit, the detailed model of BM4 is not presented here.

To quantify the performances of the different controllers, the following indices are defined:

\emph{1) Total Loss of Load (in MWh):}
\begin{align}
  \mathrm{LoL} = \int_{0}^{t_\mathrm{f}} P_\mathrm{D}^\mathrm{sh}(t) dt
\end{align}

\emph{2) Root-Mean-Square Frequency Deviation (in Hz):}
\begin{align}
  \mathrm{FD} = \sqrt{\frac{1}{t_\mathrm{f}}\int_{0}^{t_\mathrm{f}} \omega(t)^2 dt}
\end{align}

\begin{table}[tb]\scriptsize\color{black}
	\processtable{\color{black}Comparison of the control performance indices of the benchmarking and proposed methods \label{tab:perf}}
	{\begin{tabular}{@{\extracolsep{\fill}}llll}
		\toprule
		 Method   & Loss of Load  & Frequency Deviation    & Water Stage  Feasibility   \\
                  &   (MWh)       &  (Hz)                  &                            \\
		\midrule
		BM1       &  $-$          & 0.0027                 & \emph{Infeasible}    \\
		BM2       &  \emph{78}    & 0.0071                 & Feasible             \\
		BM3       &  58.83        & \emph{0.1003}          & Feasible             \\
		BM4       &  36           & 0.0887                 & Feasible             \\
		Proposed  &  \emph{32}    & 0.0867          & Feasible             \\
		\botrule
	\end{tabular}}{}
\end{table}

\begin{figure}[tb]
    {\centering
    \subfloat[]{\includegraphics[width=2.3in]{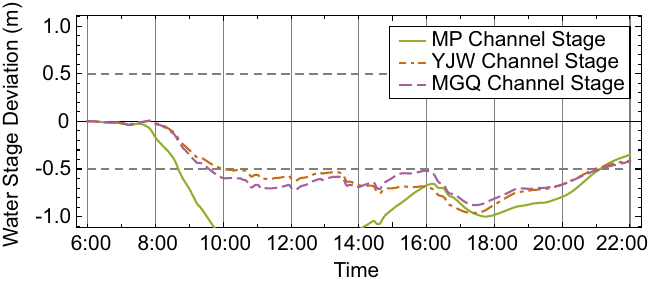}} \\
    \subfloat[]{\includegraphics[width=2.3in]{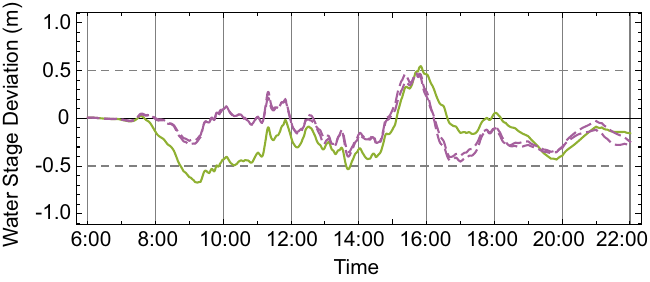}} \\
    \subfloat[]{\includegraphics[width=2.3in]{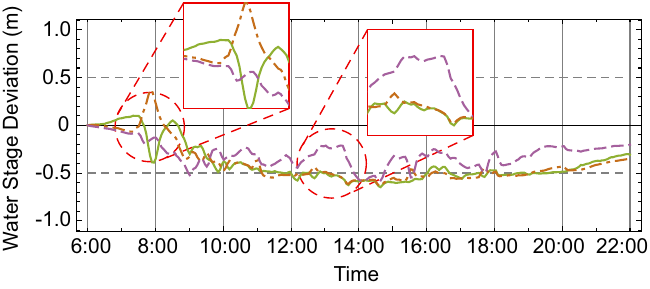}} \\
    \subfloat[]{\includegraphics[width=2.3in]{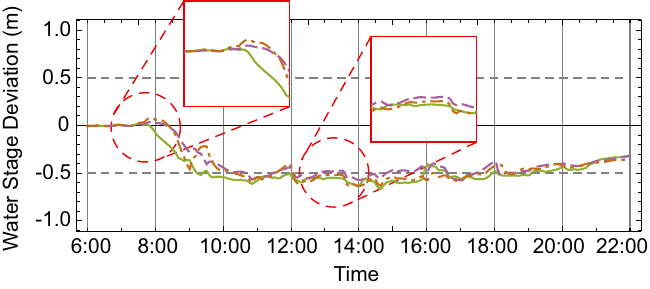}} \\}
    \caption{\color{black}Water stage deviation at the cascaded hydropower plants under different control methods}
    \figfooter{a}{\color{black}Benchmarking Method 1 (BM1)}
    \figfooter{b}{\color{black}Benchmarking Method 3 (BM3)}
    \figfooter{c}{\color{black}Benchmarking Method 4 (BM4)}
    \figfooter{d}{\color{black}The proposed control method}
	\label{fig:stage}
\end{figure}

By simulation, the performance indices of the benchmarking methods and the feasibility of the river stage constraints are listed in Table \ref{tab:perf}. In detail, because load shedding is not considered in BM1 and the energy use in stabilising the frequency is not appropriately allocated to the three cascaded plants, water used for power generation significantly excesses the river upstream inflow at the MP plant, and as the result the water stage at the MP plant descended to an unacceptably low value ($<-1.0$ m), as shown in Fig. \ref{fig:stage}(a). This violates the river operation constraints and may exhaust the water storage; therefore, it is strictly forbidden in operation.

The load shedding control (Fig. \ref{fig:logic}) in BM2 ensures that the water stage limits are not violated. However, from the water stage curves in Fig. \ref{fig:stage}(b), we can see that the cascaded plants are not coordinated at all. This causes the adjustability of the cascaded hydropower not fully exploited, leading to the largest load shedding, i.e., {\color{black}$78$} MWh, as shown in Table \ref{tab:perf}. Setting the frequency reference to the lower limit, i.e., $-0.1$ Hz, to reduce the load demand in BM3, load shedding still reaches {\color{black}$58.83$} MWh.

BM4 coordinates the cascaded plants by scheduling the base power references every $10$ minutes. Thus, the adjustability is significantly improved, indicated by the decreased load shedding shown in Table \ref{tab:perf}. However, as shown in Fig. \ref{fig:pgenschd}(a), the power outputs of the plants deviate from the scheduling within each $10$-minute interval due to the solar power and load volatility and the consequent control actions of the AGC. This portion of the power adjustment is not coordinated in BM4, in contrast to the proposed method, again revealed by the differences between the three water stage curves of different plants, as plotted in Fig. \ref{fig:stage}(c). In other words, the adjustability of the hydropower can be further exploited.

\subsection{Simulation Result of the Proposed Regulation Method}
\label{sec:result}

\begin{figure}[tb]
    {\centering
    \subfloat[]{\includegraphics[width=2.3in]{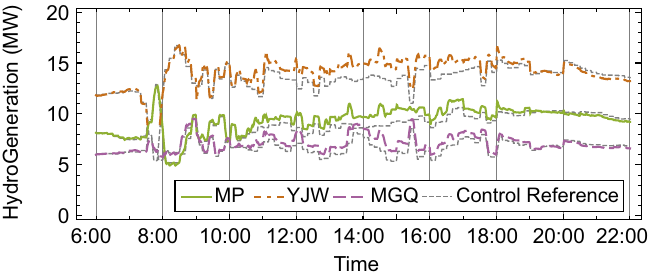}} \\
    \subfloat[]{\includegraphics[width=2.3in]{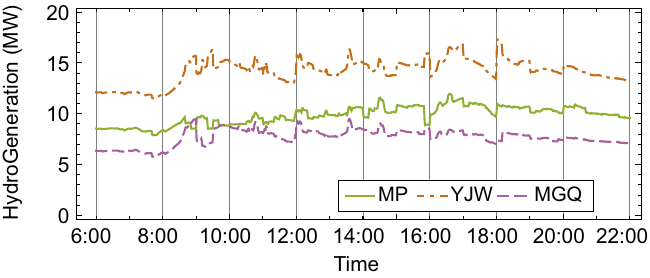}} \\}
    \caption{\color{black}Power generation of the cascaded hydropower plants under BM4 and the proposed tertiary regulation method}
    \figfooter{a}{\color{black}BM4, where the references given by scheduling are shown as dashed gray curves}
    \figfooter{b}{\color{black}The proposed method}
	\label{fig:pgenschd}
\end{figure}

\begin{figure}[tb]
    {\centering
	\includegraphics[width=2.3in]{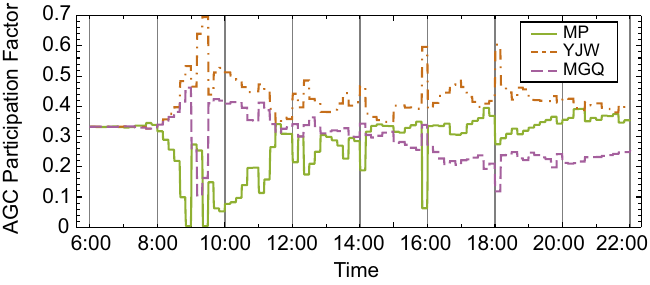}\\}
	\caption{\color{black}AGC participation factors obtained by the proposed method}
	\label{fig:weight}
\end{figure}

\begin{figure}[tb]
    {\centering
	\includegraphics[width=2.3in]{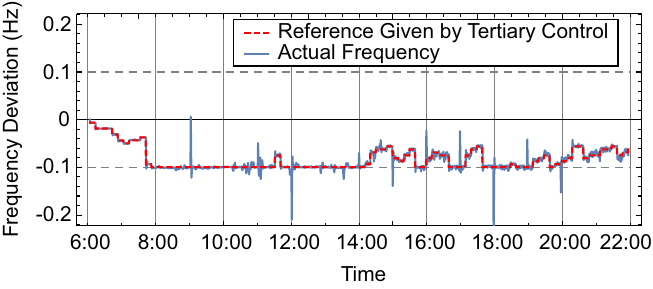}\\}
	\caption{\color{black}System frequency deviation under the proposed method}
	\label{fig:freqmpc}
\end{figure}

\begin{figure}[tb]
    {\centering
	\includegraphics[width=2.3in]{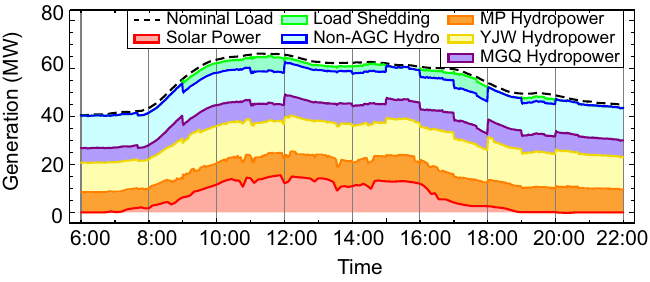}\\}
	\caption{\color{black}Total power generation under the proposed method with different components plotted in different colors}
	\label{fig:tot}
\end{figure}

By simulation of the proposed method, the AGC participation factors, as the outputs of the proposed tertiary regulation method, are plotted in Fig. \ref{fig:weight}. The system frequency reference and the actual deviation are shown in Fig. \ref{fig:freqmpc}. The power generations of the plants are shown in Fig. \ref{fig:pgenschd}(b), and a comparison of the total power generation and load shedding to the load demand is given in Fig. \ref{fig:tot}. The water stages of the plants are given in Fig. \ref{fig:stage}(d).

\begin{figure}[tb]
	{\centering
	\includegraphics[width=2.4in]{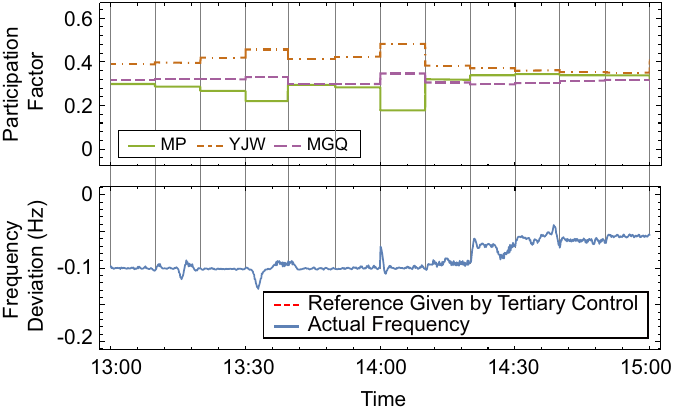}\\}
	\caption{\color{black}Events of changing AGC participating factors and system frequency deviation between 13:00 and 15:00 under the proposed tertiary control}
	\label{fig:freq2hour}
\end{figure}

The overall performance indices of the proposed method are listed in Table \ref{tab:perf} compared to those of the benchmarking methods. As can be seen, the proposed method further reduces the total load loss and frequency deviation compared all the benchmarking methods. {\color{black} The maximal transient frequency deviation over the whole day is $+0.006$ Hz and $-0.262$ Hz, which satisfies the system's operational margin, i.e., $\pm0.5$ Hz.}

{\color{black}
If we do not apply the pv smoothing curtailment control but taking the risk of instability caused by pv volatility, the total load shedding is $17$ MWh. The result with smoothing curtailment control, i.e., $32$ MWh, is to some extent larger. Alternatively, if we simply disconnect the pv plant from the grid, due to the power supply shortage, load shedding jumps to $148$ MWh, which is very socioeconomically costly. This result proves that tracking the pv prediction curve by smoothing curtailment is a sweet-point tradeoff between frequency stability and power supply reliability.
}

The following phenomena exemplify how the proposed method coordinates the cascaded plants. From Figs. \ref{fig:weight} and \ref{fig:stage}(d), at 10:00 am, the stage at the MP plant descends to the lower limit, causing the plant to lose its power adjustment ability. This is caused by the drop in the upstream inflow shown in Fig. \ref{fig:scen}(a). In this situation, the proposed controller lowers the AGC participation factor of MP to almost zero in response. Meanwhile, the other two plants take the responsibility of mitigating solar power volatility in the secondary control, as their AGC participation factors are raised. After 20:00 since solar power drops to zero and does not fluctuate, the proposed controller slowly adjusts the participation factors to drive the water stages of the plants to rise to the nominal value slowly and synchronously.

Then, comparing Figs. \ref{fig:stage}(c) and \ref{fig:stage}(d), especially the magnified parts, the difference between the three water stage curves of the cascaded plants becomes smaller than that under BM4. This outcome reveals that the proposed method better coordinates the cascaded plants. As a result, the load loss and frequency deviation are reduced, as compared in Table \ref{tab:perf}.

Moreover, comparing Figs. \ref{fig:pgenschd}(a) and \ref{fig:pgenschd}(b), we can see that the trajectories of hydropower generation under the proposed method are smoother than those of BM4. This result again shows that the proposed method offers a better coordination of the cascaded plants over short time periods. Obviously, this improvement is achieved by jointly considering the multi-timescale dynamics.

Finally, to analyse the impact of changing AGC participating factors on system frequency, we magnify the factors and frequency deviation between 13:00 and 15:00 in Fig. 14 with the events of changing AGC participating factors labeled. We can see the frequency fluctuates around the reference smoothly, and no observable disruption is observed at the events of changing the participating factors. This can be attributed to the fact that the frequency a power system is mainly stabilised by the primary control, and the AGC is designed to removes all stationary-state frequency deviation of primary control as long as the adjusting abilities of the hydro plants are preserved, which is guaranteed by the inequality constraints in the proposed tertiary controller.

\subsection{Numerical Discussions}
\label{sec:discussion}

\subsubsection{Impact of objective weights on control performance}
\label{sec:weights}

Section \ref{sec:objective} briefly discusses on how to choose the weights in the objective function (\ref{eq:obj}) of the proposed tertiary controller. Here we numerically demonstrate how these weights affect the control performance. With the weights perturbed, the corresponding components of the objective and performance indices are summarised in Table \ref{tab:weights}. As we can see, with one weight increased, the corresponding objective component decreases but others increase. The result shows that in engineering practice, the tradeoff between different objective components can be flexible adjusted by these weights.

Note that despite given different objective weight settings, the frequency stability of the system is not impacted. This is attributed to the fact that frequency stability is guaranteed by the primary and secondary control. With inequality constraints, the proposed tertiary control maintains the primary and secondary frequency control ability for in the longer time scale, which is not affected by choosing different objective weight factors.

\begin{table}[tb] \color{black}
	\processtable{ \color{black}Impact of objective weights on the control performance \label{tab:weights}}
	{\begin{tabular}{@{\extracolsep{\fill}}llllll}
		\toprule
		Weights in obj. (\ref{eq:obj})                            & RMS obj components                       & Load loss  & $\Delta$Freq    &   \\
        $\{\lambda_1,\lambda_2,\lambda_3,\lambda_4,\lambda_5\}$   & $\{J_1,J_2,J_3,J_4,J_5\}$  &   (MWh)        &(Hz)                           &                            \\
		\midrule
		\{10,10,1,1000,10\}*      &  \{0.0082,2.75,1.41,0,0.367\}             & 32              & 0.087           \\
		\{{\bf40},10,1,1000,10\}  &  \{{\bf0.0034},3.40,1.54,0,0.364\}        & 39              & {\bf0.048}         \\
		\{10,{\bf40},1,1000,10\}  &  \{0.0083,{\bf2.75},1.41,0,0.367\}        & {\bf31}         & 0.088      \\
        \{10,10,{\bf5},1000,10\}  &  \{0.0098,2.75,{\bf1.28},0,0.375\}        & 31              & 0.099         \\
        \{10,10,5,1000,{\bf50}\}  &  \{0.0091,2.75,1.45,0,{\bf0.360}\}        & 31              & 0.093         \\
		\botrule
	\end{tabular}}{
    * the original weights used in Section \ref{sec:result}.}
\end{table}

\subsubsection{Robustness of the proposed controller against perturbed modeling parameters}
\label{sec:modelerror}

\begin{figure}[tb]
	\centering
	\includegraphics[width=3.3in]{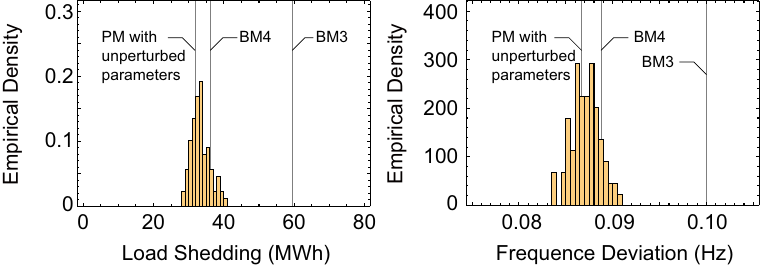}
	\caption{\color{black}Empirical densities of the overall load shedding and RMS frequency deviation under the proposed control method with $100$ sets of randomly perturbed modeling parameters}
	\label{fig:distpara}
\end{figure}

To test the robustness of the proposed control method against inaccurate system model, we independently perturb the modeling parameters in the MPC (including river section geometry and hydropower production efficiency) by $\pm5\%$ with uniform distributions. Overall $100$ sets of perturbed parameters are created and used for simulation. Except for the parameters in MPC formulation, system parameters for simulation remain untouched. Empirical densities of the obtained performance indices of load shedding and RMS frequency deviation are shown in Fig. \ref{fig:distpara}. Despite significantly disturbed parameters, the performance of the proposed control method remains relatively consistent. This can also be attributed to the inherent robustness of the MPC.

\subsection{Simulation Results of the Proposed Regulation Method under Various Scenarios}
\label{sec:samples}

\begin{figure}[tb]
    {\centering
    \subfloat[]{\includegraphics[width=2.3in]{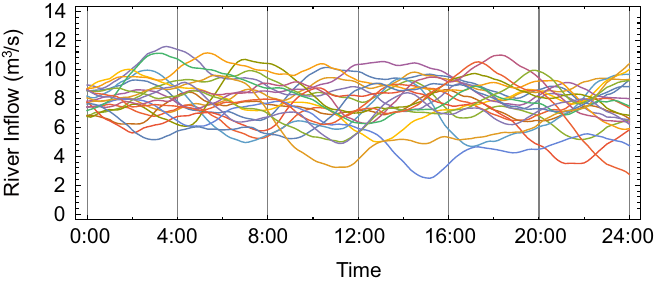}}  \\
    \subfloat[]{\includegraphics[width=2.3in]{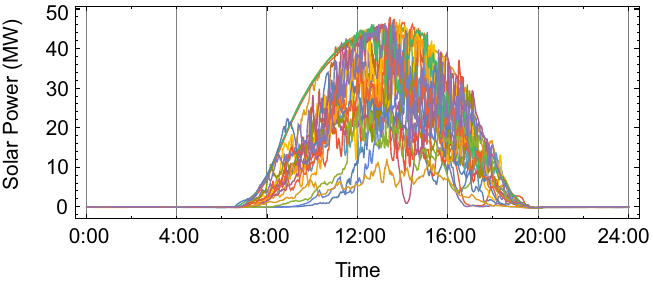}}      \\
    \subfloat[]{\includegraphics[width=2.3in]{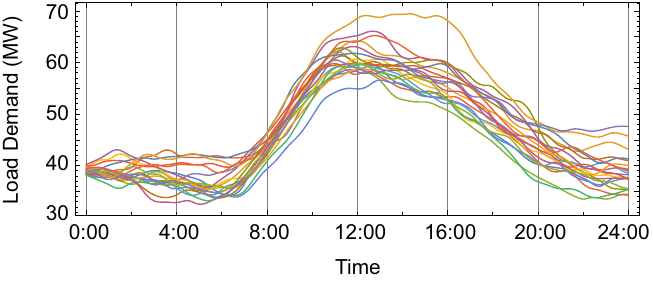}}    \\ }
    \caption{Twenty out of one hundred simulation scenarios used in the case study in Section \ref{sec:samples} }
    \figfooter{a}{River upstream inflow}
    \figfooter{b}{{\color{black}Available} solar power generation}
    \figfooter{c}{Load demand}
	\label{fig:many}
\end{figure}

To validate the improvement in the proposed method more comprehensively, especially considering various operational conditions in a real-life system, we tested it under $100$ different scenarios. Each scenario represents one day from 6:00 to 22:00. The solar power curves used are empirical data recorded in the Xiaojin system from April to July, 2018. For visualisation, $20$ out of the $100$ scenarios used for simulation are plotted in Fig. \ref{fig:many}. Specifically, the data of {\color{black}available pv power} shown in Fig. \ref{fig:many}(b) are recorded on $100$ consecutive days from January to April, 2018. {\color{black}Smoothing control in Section \ref{sec:smoothing} is applied in simulation.}

\begin{figure}[tb]
	{\centering
	\includegraphics[width=2.3in]{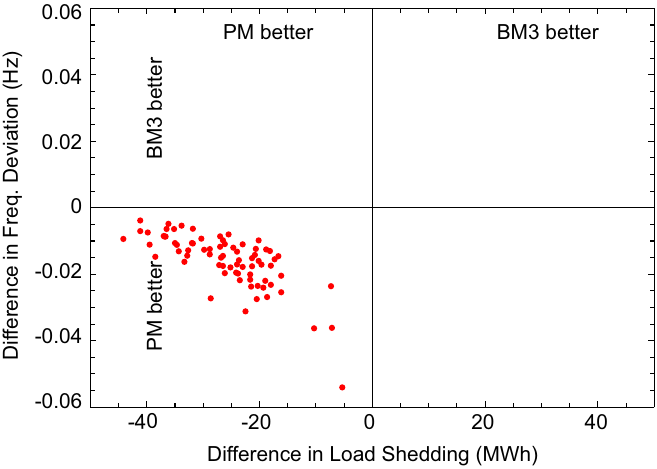}\\}
	\caption{\color{black}Difference in the performance indices between the proposed method (PM) and BM3 under $100$ scenarios. Each mark represents one scenario}
	\label{fig:samplesbm3}
\end{figure}

\begin{figure}[tb]
	{\centering
	\includegraphics[width=2.3in]{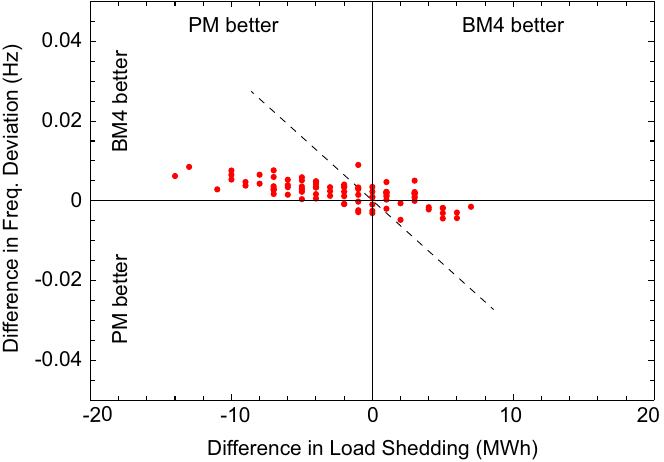}\\}
	\caption{\color{black}Difference in the performance indices between the proposed method (PM) and BM4 under $100$ scenarios. Each mark represents one scenario} 
	\label{fig:samples}
\end{figure}

Simulation results of the proposed method are compared to those of BM3 and BM4.
The differences between the total load loss and frequency deviation of the proposed method and those of BM3 and BM4 are respectively plotted as one mark per scenario in Figs.  \ref{fig:samplesbm3} and \ref{fig:samples}. Seen from Fig. \ref{fig:samplesbm3}, compared to BM3 the proposed controller performers dramatically better in terms of reducing total load loss. Comparing the proposed method and BM4, as observed in Fig. \ref{fig:samples}, most points appear on the left halfplane, meaning that in most scenarios the proposed method still outperforms BM4 in terms of load loss while keeping frequency deviation almost constant.

Based on all the results  above, we can conclude that the proposed regulation method that considers multi-timescale dynamics better coordinates the cascaded plants and exploits the overall power adjustability compared to the benchmarking methods, and thus improves the power supply ability of islanded power systems with cascaded run-of-the-river hydropower.

\section{Conclusions}

A coordinated tertiary regulation approach for islanded power systems with cascaded run-of-the-river hydropower and volatile generations is proposed. An MPC is established to dynamically adjust the AGC participation factors to coordinate the cascaded plants. A simulation of a real-life system shows that the proposed regulation method significantly reduces load loss compared to that yielded by the other regulation methods.

Currently, the proposed method has not taken the statistical characteristics of solar and wind power into modeling and optimisation, but recent work has shown that considering stochastic characteristics improves the performance \cite{qiu2020stochastic}. In future studies, taking the stochastic characteristics of the volatile generations and load into the tertiary regulation could be a promising work direction.

\section{Acknowledgments}

Financial support came from the National Key Research \& Development Program of China (2018YFB0905200), the National Natural Science Foundation of China (51907099, 51677100, 51761135015), and the China Postdoctoral Science Foundation (2019M650676).

\section{Statements}

\subsection{Conflict of Interest}

There is no conflict of interest to report.

\subsection{Data Availability }
The data that support the findings of this study are available from the corresponding author upon reasonable request.

\bibliographystyle{iet}
\bibliography{HYMGAGC}

\section{Appendix A: Polynomial Approximation of the Power System Dynamics}
\label{sec:polyapp}

In Section \ref{sec:modelelectric} polynomial functions (\ref{eq:polyapp}) are used to approximate the mean hydropower outputs over the time interval of $t\in[0,T]$, which are exemplified in Table \ref{tab:poly}. This appendix elaborates on how to obtain the polynomial surrogate models.

Because Legendre polynomials are orthogonal over a finite input interval \cite{Xiu2010,Qiu2020Explicit}, leading to the optimal approximation in terms of $L_2$-norm, this work choose Legendre polynomials as the basis in approximation.

The $j$th-order Legendre polynomial $P_j(x)$ for $x$ is defined as
\begin{align}
    P_j(x) = \frac{1}{2^j} \sum_{k=0}^{j}\binom{j}{k}^2 (x-1)^{j-k} (x+1)^{k}, \label{eq:legendre}
\end{align}
\noindent
which admit the orthogonality over interval $\mathbb{I}=[-1,1]$, as
\begin{align}
    \langle \phi_{i}(x), \phi_{j}(x) \rangle_{L^2(\mathbb{I})} :=  \int_{\mathbb{I}} \phi_{i}(x) \phi_{j}(x) dx =  \frac{2}{2i+1} \delta_{i j},  \label{eq:motho}
\end{align}
\noindent
where $\langle\cdot,\cdot\rangle_{L^2(\mathbb{I})}$ represents the inner product over $\mathbb{I}$; $\delta_{jk}$ is the Kronecker delta function.

Specifically, the Legendre polynomials up to order $2$ are
\begin{align}
    P_0(x) = 1,\ P_1(x)=x,\ P_2(x)=-0.5+1.5^2. \nonumber
\end{align}

In finding the surrogate models for the power system dynamics in Section \ref{sec:modelelectric}, the function to be approximated, i.e., $\bar{P}_{\mathrm{G}i}(\cdot)$, has multidimensional inputs, denoted as vector $\bm{q}$,
\begin{align}
  \bm{q} := \left[ \bar{P}_{\mathrm{S}}, \bar{P}_{\mathrm{D}} , \bm{p}^\mathrm{T} \right]^\mathrm{T} = \left[ \bar{P}_{\mathrm{S}}, \bar{P}_{\mathrm{D}}, \bm{c}_{\mathrm{G}}^{\mathrm{T}}, \omega^{\mathrm{ref}}, P_\mathrm{D}^\mathrm{sh} \right] ^\mathrm{T}.
\end{align}

We first rescale all inputs into interval $\mathbb{I}=[-1,1]$, and denote the dimension of $\bm{q}$ as $m$. Given the approximation order $N_i$ for each scaler input $q_i$, the multivariate Legendre polynomial basis $\left\{ \phi_{\bm{j}}(\bm{q}) \right\}$ for the vector input $\bm{q}$ is constructed as the following tensor product,

\noindent
\begin{align}
    \hspace{-2pt}\left\{ \phi_{\bm{j}}(\bm{q}) \right\} = \left\{ P_{j_1}(q_1) P_{j_2}(q_2) \cdots P_{j_m}(q_{j_{m}}) : j_i \le N_i \right\},  \hspace{-2pt} \label{eq:wick}
\end{align}

\noindent
where $\bm{j} \triangleq \left[ j_1,j_2,\ldots,j_{M} \right]^{\mathrm{T}}$ is the multi-dimensional index.

Then, the function $\bar{P}_{\mathrm{G}i}(\bm{q})$ is approximated as the following polynomial expansion, denoted $\bar{P}^*_{\mathrm{G}i}(\bm{q})$, as
\begin{align}
   \bar{P}_{\mathrm{G}i}(\bm{q}) \approx \bar{P}^*_{\mathrm{G}i}(\bm{q}) := \sum_{\bm{j}} f_{\bm{j}} \phi_{\bm{j}}(\bm{q}), \label{eq:pce}
\end{align}
\noindent
where $f_{\bm{j}}$ are coefficients to be calculated.

Generally, the Galerkin method (GM) or the collocation method (CM) can be used to find the coefficients of a polynomial approximation \cite{Xiu2010,Hockenberry2004,Qiu2020Explicit}. Because the CM can treat the function to be approximated, i.e., $\bar{P}_{\mathrm{G}i}(\bm{q})$, as a black box in a simulator such as the PSS/E whereas the GM relies on a white-box model with detailed mathematical formulation, we adopt the CM to find the coefficients in the polynomial surrogate model in this work.

First, define the collocation point set, as know as the Gaussian quadrature point set \cite{Xiu2010,Hockenberry2004}, to contain all the zeros of the product of the $(N_i+1)$th Legendre polynomials of the inputs in $\bm{q}$, as
\begin{align}
    \big\{ \hat{\bm{q}}_i : P_{N_1+1}(q_1) P_{N_2+1}(q_2) \cdots P_{N_{m}+1}(q_m) = 0 \big\}_{i=1}^{N_b}, \label{eq:cps}
\end{align}
\noindent
where $N_b$ is the size of collocation point set, which equals the size of the basis (\ref{eq:wick}).

For each collocation point $\hat{\bm{q}}_i$ in (\ref{eq:cps}), using (\ref{eq:loadnorm}) to create the corresponding trajectories of solar power and load demand and set the tertiary control variables according to $\hat{\bm{q}}_i$. Then, using power system simulation software to find the values of $P_{\mathrm{G}i}(\hat{\bm{q}}_i)$.
After the simulations of all collocation points are finished, the coefficients in the polynomial approximation model (\ref{eq:pce}) are computed by
\begin{align}
    \left[  f_1,  \ldots, f_{N_b}  \right]^{\mathrm{T}} =
    \bm{A}^{-1}
    \left[  P_{\mathrm{G}i}(\hat{\bm{q}}_1) , \ldots, P_{\mathrm{G}i}(\hat{\bm{q}}_{N_b}) \right]^{\mathrm{T}}, \label{eq:cm}
\end{align}
\noindent
where $\bm{A}$ is a constant matrix determined by the basis (reindexed as $\left\{ \phi_i(\cdot) \right\}_{i=1}^{N_b}$) and the collocation points (\ref{eq:cps}), as
\begin{align}
    \bm{A} =
    \begin{bmatrix}
        \phi_1(\hat{\bm{q}}_1)     & \cdots & \phi_{N_b}(\hat{\bm{q}}_1) \\
        \vdots                       & \ddots & \vdots   \\
        \phi_1(\hat{\bm{q}}_{N_b}) & \cdots &  \phi_{N_b}(\hat{\bm{q}}_{N_b})
         \\
    \end{bmatrix}.  \label{eq:cmat}
\end{align}

Alternatively, due to that the high-order cross terms in the polynomial expression generally have less impact on the precision of the approximation \cite{Xiu2010}, we can use the \emph{Smolyak adaptive sparse algorithm}\cite{CONRAD2013ADAPTIVE} to incrementally add higher-order cross terms in an adaptive way instead of adding all high-order cross terms at once. The precision of the approximation in terms of $L_2$-norm is monitored by the iterative changes in the coefficients, and only the terms with a noticeable impact are added. This process terminates when the preset error tolerance $\epsilon$ is reached. Interested readers are refereed to the detailed theory of adaptive polynomial approximation in \cite{CONRAD2013ADAPTIVE} as well as its applications in \cite{Qiu2020nonintrusive}.

\section{Appendix B: Test System Parameters}
\label{sec:parameters}

\begin{table}[h]
	\processtable{River section data in the cascaded hydropower system \label{tab:river}}
	{\begin{tabular}{@{\extracolsep{\fill}}llllll}
		\toprule
		\# & Section Type      & Length (m)    & Width (m)     & Slope (\%)    & Friction Factor  \\
		\midrule
		1   &  Natural River   & 15000   & 14.34  & 1.35  & 0.030      \\
		2   &  Natural River   & 10000   & 18.23  & 1.23  & 0.030      \\
		3   &  Channel         & 10000   & 3.30   & 0.07  & 0.012      \\
		4   &  Natural River   & 1800    & 20.50  & 1.25  & 0.030      \\
		5   &  Natural River   & 11000   & 20.84  & 1.03  & 0.030      \\
		6   &  Channel         & 11400   & 3.30   & 0.07  & 0.012      \\
		7   &  Natural River   & 11000   & 24.15  & 1.98  & 0.030      \\
		8   &  Natural River   & 8000    & 21.64  & 2.01  & 0.030      \\
		9   &  Channel         & 8300    & 3.30   & 0.07  & 0.012      \\
		10  &  Natural River   & 7000    & 19.49  & 1.55  & 0.030      \\
		\botrule
	\end{tabular}}{}
\end{table}

\begin{table}[h]
	\processtable{Brief data of the hydropower plants in the simulation \label{tab:plant}}
	{\begin{tabular}{@{\extracolsep{\fill}}llllll}
		\toprule
		$\#$ & Name  & Rated Power (MW)  & Initial Power (MW)     & Initial Head (m)  \\
		\midrule
		1   &  MP    & 15   & 8.47       & 125.6       \\
		2   &  YJW   & 20    & 12.24      & 181.3       \\
		3   &  MGQ   & 12    & 6.14       & 91.0       \\
		\botrule
	\end{tabular}}{}
\end{table}

\begin{table}[h]
	\processtable{Key parameters of the HYGOVM turbine-governor models\label{tab:turbine}}
	{\begin{tabular}{@{\extracolsep{\fill}}llll}
		\toprule
		Parameter                                       & \phantom{$-$}\#1 MP  & \phantom{$-$}\#2 YJW & \phantom{$-$}\#3 MGQ   \\
		\midrule
		Rated power ($\text{MW}$)                       & \phantom{$-$}15      & \phantom{$-$}20      & \phantom{$-$}12        \\
        Rated discharge ($\text{m}^3/\text{s}$)         & \phantom{$-$}11      & \phantom{$-$}12.8    & \phantom{$-$}13.3    \\
        Rated head ($\text{m}$)                         & \phantom{$-$}135.0   & \phantom{$-$}184.0   & \phantom{$-$}91.0    \\
        No load flow ($\text{p.u.}$)                    & \phantom{$-$}0.08    & \phantom{$-$}0.08    & \phantom{$-$}0.08    \\
        Permanent droop                                 & \phantom{$-$}0.05    & \phantom{$-$}0.05    & \phantom{$-$}0.05    \\
        Temporary droop                                 & \phantom{$-$}0.10    & \phantom{$-$}0.10    & \phantom{$-$}0.10    \\
        \color{black} Deadband                           & \color{black}$\pm0.05$ Hz & \color{black}$\pm0.05$ Hz  & \color{black}$\pm0.05$ Hz \\
        Governor time constant ($\text{s}$)             & \phantom{$-$}5.0     & \phantom{$-$}5.0     & \phantom{$-$}5.0     \\
        Filter time constant ($\text{s}$)               & \phantom{$-$}0.05    & \phantom{$-$}0.05    & \phantom{$-$}0.05    \\
        Servo time constant ($\text{s}$)                & \phantom{$-$}0.5     & \phantom{$-$}0.5     & \phantom{$-$}0.5    \\
        Max gate opening rate ($\text{s}^{-1}$)         & \phantom{$-$}0.1     & \phantom{$-$}0.1     & \phantom{$-$}0.1     \\
        Max gate closing rate ($\text{s}^{-1}$)         & $-$0.125  & $-$0.125  & $-$0.125  \\
        Max buffer gate opening rate ($\text{s}^{-1}$)  & \phantom{$-$}0.1     & \phantom{$-$}0.1     & \phantom{$-$}0.1     \\
        Max buffer gate closing rate ($\text{s}^{-1}$)  & $-$0.05   & $-$0.05   & $-$0.05   \\
        Penstock length/cross section ($\text{m}^{-1}$) & \phantom{$-$}39.23   & \phantom{$-$}31.99   & \phantom{$-$}26.47   \\
        Penstock head loss factor ($\text{s}^{2}/\text{m}^5$) & \phantom{$-$}0.00043   & \phantom{$-$}0.00022   & \phantom{$-$}0.00077   \\
		\botrule
	\end{tabular}}{}
\end{table}

\begin{table}[h]\color{black}
	\processtable{\color{black}Per unit network parameter\label{tab:netpara}}
	{\begin{tabular}{@{\extracolsep{\fill}}lllll}
		\toprule
		Bus From  \hspace{6pt}    & Bus To \hspace{12pt} &R   \hspace{32pt}    & X  \hspace{32pt} & B \\
		\midrule
		1	&2	      &0.0025	&0.0600	&0 \\
        1	&59	      &0.0005	&0.0037	&0.0074\\
        1	&62	      &0.0020	&0.0106	&0.0249\\
        2	&3	      &0.0007	&0.0021	&0\\
        2	&6	      &0.0006	&0.0019	&0\\
        3	&4	      &0.0290	&0.641	&0\\
        3	&26	      &0.0236	&0.0705	&0\\
        6	&7	      &0.0100	&0.240	&0\\
        9	&101	  &0.264	&0.371	&0\\
        19	&101	  &0.0044	&0.0062	&0\\
        21	&22	      &0.0794 	&1.91	&0\\
        21	&102	  &0.209	&0.294	&0\\
        21	&104	  &0.0336	&0.101	&0\\
        23	&104	  &0.026	&0.0776	&0\\
        24	&25	      &0.0794 	&1.91	&0\\
        24	&28	      &0.213	&0.297	&0\\
        26	&27	      &0.0100	&0.240	&0\\
        26	&28	      &0.0100	&0.240	&0\\
        27	&33	      &0.121	&0.135	&0\\
        28	&104	  &0.123	&0.367	&0\\
        33	&34	      &0.125 	&3.00	&0\\
        35	&36  	  &0.100	&2.40	&0\\
        35	&37	      &0.100	&2.40	&0\\
        35	&41	      &0.113	&0.126	&0\\
        35	&47	      &0.617	&0.866	&0\\
        36	&49	      &0.200	&0.100	&0\\
        37	&43	      &0.100	&0.100	&0\\
        41	&42	      &0.0625	&1.50	&0\\
        41	&104	  &0.100	&0.100	&0\\
        43	&44	      &0.250	&6.00	&0\\
        47	&48	      &0.0526	&1.26	&0\\
        49	&50	      &0.100    &2.40	&0\\
        59	&61	      &0.0200   &0.480	&0\\
        62	&63	      &0.0167	&0.400	&0\\
        101	&102	  &0.154	&0.216	&0\\
		\botrule
	\end{tabular}}{}
\end{table}

\end{document}